\documentclass[sn-mathphys]{sn-jnl}
\jyear{2022}%

\usepackage{amssymb}
\usepackage{amsmath,amsfonts,amssymb, amsthm}
\usepackage{extarrows}
\usepackage{paralist}
\usepackage{graphics}
\usepackage{epsfig}
\usepackage{psfrag}
\usepackage{subeqnarray}
\usepackage{cases}
\usepackage{bm}
\usepackage{cleveref}
\usepackage{bibentry}
\usepackage{appendix}
\usepackage{xcolor}
\usepackage{enumerate}
\usepackage{color}
\usepackage[all,cmtip,line]{xy}
\usepackage{array,setspace}

\newcommand{\dd}{{\rm d}}
\newcommand{\xx}{\mathbf{x}}
\newcommand \RR{{\mathbb{R}}}
\newcommand{\nnu}{{\boldsymbol\nu}}
\newcommand{\uu}{\mathbf{u}}
\newcommand{\vv}{\mathbf{v}}
\newcommand \der{\partial}
\newcommand \vphi{\varphi}

\newcommand \gam{\gamma}

\newcommand{\dr}{{\rm d}}
\newcommand{\D}{{\rm D}}

\newcommand{\bR}{ {\mathbb{R}}} 

\newcommand{\defd}{:=}

\newcommand{\bea}{{\begin{eqnarray}}}
\newcommand{\eea}{{\end{eqnarray}}}

\newcommand{\FF}{{\boldsymbol F}}

\newcommand \ol{\overline}

\newcommand \Sonic{\Gamma_{\rm sonic}}
\newcommand \Shock{\Gamma_{\rm shock}}
\newcommand \Wedge{\Gamma_{\rm wedge}}

\newcommand \shock{\Gamma_{\rm shock}}

\newcommand \ivphi{\varphi_0}

\newcommand \irho{\rho_0}

\newcommand{\pSi}{\varphi}
\newcommand{\PtIncW}{{P_0}}
\newcommand{\PtUpL}{{P_1}}
\newcommand{\PtLwL}{{P_2}}
\newcommand{\PtLwR}{{P_3}}
\newcommand{\PtUpR}{{P_4}}
\newcommand \bPhi{\bar{\Phi}}
\newcommand{\xxi}{{\boldsymbol\xi}}

\newcommand{\Gso}{\Gamma_{\text{\rm sonic}}}

\newcommand{\Gsh}{\Gamma_{\text{\rm shock}}}

\newcommand{\bn}{{\boldsymbol\nu}}

\newcommand{\mR}{\mathbb{R}}

\newcommand{\cxi}{{\xi_1}}
\newcommand{\ceta}{{\xi_2}}

\renewcommand{\r}{\rho}

\newcommand \Symm{\Gamma_{\rm sym}}

\numberwithin{equation}{section}
\numberwithin{figure}{section}
\theoremstyle{plain}
\newtheorem{definition}{Definition}[section]
\newtheorem{theorem}{Theorem}[section]

\newtheorem{problem}{Problem}[section]

\begin{document}

\title[Two-Dimensional Riemann Problems]{Two-Dimensional Riemann Problems: \\ Transonic Shock Waves\\ and Free Boundary Problems}


\author{\fnm{Gui-Qiang G.} \sur{Chen}}\email{chengq@maths.ox.ac.uk}


\affil{\orgdiv{Oxford Centre for Nonlinear PDE, Mathematical Institute}, \\ \orgname{University of Oxford}, \orgaddress{\city{Oxford}, \postcode{OX2 6GG}, \country{UK}}}


\abstract{We are concerned with global solutions of multidimensional Riemann problems
for nonlinear hyperbolic systems of conservation laws, focusing on their global configurations and structures.
We present some recent developments in the rigorous analysis of two-dimensional Riemann problems
involving transonic shock waves through several prototypes of hyperbolic systems of conservation laws
and discuss some further multidimensional Riemann problems and related problems for nonlinear
partial differential equations (PDEs).
In particular, we present four different two-dimensional Riemann problems
through these prototypes of hyperbolic systems
and show how these Riemann problems can be reformulated/solved as
free boundary problems with transonic shock waves as free boundaries
for the corresponding nonlinear conservation laws of mixed elliptic-hyperbolic type
and related nonlinear PDEs.
}


\maketitle

\section{Introduction}\label{sec1}

We are concerned with global solutions of multidimensional (M-D) Riemann problems
for nonlinear hyperbolic systems of conservation laws, focusing on
their global configurations and structures.
In this paper, we present some recent developments in the
rigorous analysis of two-dimensional (2-D) Riemann problems
involving transonic shock waves (shocks, for short)
through several prototypes of hyperbolic systems of conservation laws
and discuss some further M-D Riemann problems and related problems for nonlinear
partial differential equations (PDEs).
These Riemann problems can be reformulated as free boundary problems with transonic shocks
as free boundaries
for the corresponding nonlinear conservation laws of mixed elliptic-hyperbolic type
and related nonlinear PDEs.

The study of Riemann problems has an extensive history, which dates back to the pioneering
work of Riemann \cite{Ri} in 1860.
For the one-dimensional (1-D) Riemann problem, a theory has been established for the appropriate amplitude
of the Riemann data for general strictly hyperbolic systems ({\it cf}. \cite{Lax,Liu1})
and for general Riemann data
for the compressible Euler equations ({\it cf}. \cite{CH,MP,Sm,We} and the references cited therein).
The 1-D Riemann problem has been essential in the development of the 1-D mathematical theory of hyperbolic conservation laws
and associated shock capturing methods
for the construction and computation of global entropy solutions;
see \cite{Da,Gli1,GlimmMajda,Lax1,Lax,LeV,Liu1,Shu} and the references cited therein.
More importantly, general global entropy solutions
can be locally approximated by the Riemann solutions that are regarded as
fundamental building blocks of the entropy solutions ({\it cf.} \cite{Da,Gli1,Lax,Sm}).
Moreover, the Riemann solutions usually determine the large-time asymptotic behaviors and global attractors of
general entropy solutions of the Cauchy problem.
On the other hand, it is the simplest Cauchy problem (initial value problem) whose solutions have fine explicit structures.

The M-D Riemann problems are more challenging mathematically,
and the corresponding M-D Riemann solutions are of much richer global configurations and structures;
see \cite{CH,CCY1,CCY2,CCY3,CFr,Da,GlimmK,GlimmMajda,LaxLiu,SCG,ZZ} and the references cited therein.
Thus, the Riemann solutions often serve as standard test models for analytical and numerical methods
for solving nonlinear hyperbolic systems of conservation laws and related nonlinear PDEs.
Theoretical results for first-order scalar conservation laws are available in
\cite{CH,CLT,Gu,Lind,TZ,Wa,ZZ1} and the references cited therein.
During recent decades, some significant developments for the 2-D Riemann problems for first-order hyperbolic systems
and second-order hyperbolic equations
of conservation laws have been made.
Zhang-Zheng \cite{ZZ} first considered the two-dimensional
four-quadrant Riemann problem that each jump between two neighbouring quadrants projects exactly one planar fundamental wave
and predicted that there are a total of 16 genuinely different configurations of the Riemann solutions for polytropic gas.
Schulz-Rinne \cite{SR} proved that one of them is impossible.
In Chang-Chen-Yang \cite{CCY1,CCY2}, it is first observed that, when two initially parallel slip lines are present,
it makes a difference whether the vorticity waves generated have the same or opposite sign, which, along with Lax-Liu \cite{LaxLiu},
leads to the classification with a total 19 genuinely different configurations of the Riemann solutions for the compressible Euler equations
for polytropic gas, via characteristic analysis;
also see \cite{KTa,li1998two,SCG}.
On the other hand, experimental and numerical results have shown that
many new configurations may arise from other types of Riemann problems.
In particular, the angles between two discontinuities separated by sectorial regions
in the initial Riemann data and the boundaries in the lateral Riemann data play essential roles
in forming the global Riemann solution configurations, besides the strengths of jumps in
the initial Riemann data;
see \cite{BCF-14,BD,CFr,Elling,FTB,FWB,GlimmMajda,Mach,VD,Neumann0,Neumann1,Neumann2,WC}.
In this paper,
we present four different $2$-D Riemann problems involving transonic shocks through the prototypes of nonlinear hyperbolic PDEs
and demonstrate how these Riemann problems
can be reformulated and then solved rigorously as free boundary problems for
nonlinear conservation laws of mixed elliptic-hyperbolic type and related nonlinear PDEs.
A special attention has been paid to whether/how different initial or boundary setups of the Riemann problems
affect the global Riemann solution configurations.
These are achieved by developing further the nonlinear method and related ideas/techniques introduced
in Chen-Feldman \cite{ChenFeldman1,ChenFeldman,CF-book2018} for solving free boundary problems with transonic shocks as free boundaries
for nonlinear conservation laws of mixed elliptic-hyperbolic type and related nonlinear PDEs;
also see \cite{Chen2,ChenFeldman2022}.

The organization of this paper is as follows:
In Section 2, we first show how the solutions of M-D Riemann problems for hyperbolic conservation laws
can be formulated
as the self-similar solutions
for nonlinear conservation laws of mixed elliptic-hyperbolic type
and then we introduce the notion of Riemann solutions
in the self-similar coordinates in the distributional sense.
In Section 3, we present the first 2-D Riemann problem, Riemann Problem I,
involving two shocks and two vortex sheets for the pressure gradient system
and show how Riemann Problem I can be reformulated/solved as a free boundary problem with transonic shocks
as free boundaries for a second-order nonlinear conservation law of mixed elliptic-hyperbolic type
and related nonlinear PDEs.
In Section 4, we present the second 2-D Riemann problem,
Riemann Problem II -- the Lighthill problem for shock diffraction by convex cornered wedges
through the nonlinear wave equations, and show how Riemann Problem II can be solved as another free boundary problem.
Even though both the origin and form of the nonlinear wave equations are different from those of the pressure gradient system,
the same arguments for solving the Riemann problem apply for the pressure gradient system to obtain similar results
without additional analytical obstacles;
the same is true for the Riemann problem in Section 3 for the nonlinear wave equations.
In Section 5, we present the third 2-D Riemann problem, Riemann Problem III -- the Prandtl-Meyer problem for unsteady
supersonic flow onto solid wedges through the Euler equations for potential flow
and show how Riemann Problem III can be reformulated/solved as a free boundary problem
for a second-order nonlinear conservation law of mixed elliptic-hyperbolic type.
Then, in Section 6, we present the fourth 2-D Riemann problem, Riemann Problem IV -- the von Neumann
problem for shock reflection-diffraction by wedges  for the Euler equations for potential flow,
and show how Riemann Problem IV can be solved again as a free boundary problem.
We give our concluding remarks in Section 7
and discuss several further M-D Riemann problems and related problems for nonlinear PDEs.

\section{Multidimensional Riemann Problems\\ and Nonlinear Conservation Laws of \\ Mixed Elliptic-Hyperbolic Type}

In this section, we first show how the solutions of the M-D Riemann problems for nonlinear hyperbolic conservation laws
can be formulated
as the self-similar solutions for nonlinear conservation laws of mixed elliptic-hyperbolic type,
and then introduce the notion of Riemann solutions in the self-similar coordinates
in the distributional sense.

Consider both the M-D first-order quasilinear hyperbolic systems of conservation laws of the form:
\begin{equation}\label{HCL-1}
\partial_t\bm{U}+\nabla_{\xx}\cdot \bm{F}=0 \qquad\,\,\mbox{for $t\in \mathbb{R}_+=[0,\infty)$ and $\, \xx\in \mathbb{R}^n$}
\end{equation}
with
$\bm{U}\in \mathbb{R}^m$ and nonlinear mapping $\bm{F}: \mathbb{R}^m\to \mathbb{R}^m\times \mathbb{R}^n$,
and the M-D second-order quasilinear hyperbolic equations of conservation laws of the form:
\begin{equation}\label{HCL-2}
\partial_t G_0(\partial_t\Phi, \nabla_{\xx}\Phi) +\nabla_{\xx}\cdot \bm{G}(\partial_t\Phi, \nabla_\xx\Phi)=0
 \qquad\,\,\mbox{for $t\in \mathbb{R}_+$ and $\, \xx\in \mathbb{R}^n$}
\end{equation}
with
$\Phi\in \mathbb{R}$ and nonlinear mapping $(G_0, \bm{G}): \mathbb{R}^{n+1}\to \mathbb{R}\times\mathbb{R}^n$.

\smallskip
A prototype of \eqref{HCL-1} is the full Euler equations in the conservation form \eqref{HCL-1}
with
\begin{equation}\label{HCL-1e}
{\bm{U}}:=(\rho,\rho \uu,\rho E)^\top,\quad \bm{F}:=(\rho \uu,\rho\uu\otimes \uu+pI, (\rho E+p)\uu)^\top,
\end{equation}
where  $\rho>0$ is the density, $\uu\in \mathbb{R}^n$ the velocity, $p$ the pressure,
and
$
E=\frac{\lvert \uu\rvert^2}{2}+e
$
the total energy per unit mass with the internal energy $e$ given by
$e=\frac{p}{(\gamma-1)\rho}$ for the adiabatic constant $\gamma>1$ for polytropic gases.

\smallskip
A prototype of  \eqref{HCL-2} can be derived from the Euler equations for potential flow, which is
governed by the conservation law of mass and the Bernoulli law
for the density function $\rho$ and the velocity potential $\Phi$ ({\it i.e.}, $\uu=\nabla_{\xx} \Phi$):
\begin{align}\label{1-a}
\der_t\rho+ \nabla_{\xx}\cdot (\rho \nabla_{\xx}\Phi)=0,\quad\,\,
\der_t\Phi+\frac 12\lvert\nabla_{\xx}\Phi\rvert^2+h(\rho)=B,
\end{align}
where $B$ is the Bernoulli constant
and $h(\rho)$ is given by
\begin{equation}\label{1-c}
h(\rho)=\frac{\rho^{\gam-1}-1}{\gam-1} \qquad\,\, \mbox{for the adiabatic exponent $\gam>1$.}
\end{equation}
By \eqref{1-a}--\eqref{1-c},
$\rho$ can be expressed as
\begin{equation}\label{1-b1}
\rho(\der_t\Phi,\nabla_{\bf x}\Phi)=h^{-1}(B-\der_t\Phi-\frac 12\lvert \nabla_{\bf x}\Phi\rvert^2).
\end{equation}
Then system \eqref{1-a} can be rewritten as the second-order nonlinear wave equation as in
\eqref{HCL-2} with
\begin{equation}\label{HCL-2e}
(G_0, \bm{G})=(\rho(\der_t\Phi, \nabla_{\bf x}\Phi),
\rho(\der_t\Phi, \nabla_{\bf x}\Phi)\nabla_{\bf x}\Phi)
\end{equation}
and $\rho(\der_t\Phi, \nabla_{\bf x}\Phi)$ determined by \eqref{1-b1}.

\smallskip
A {\it standard Riemann problem} for \eqref{HCL-1}  is a special Cauchy problem:
\begin{equation}\label{HCL-1R}
\bm{U}\lvert_{t=0}
= \bm{U}_0(\xx)
\end{equation}
so that the initial data function $\bm{U}_0(\xx)$ is invariant under the self-similar scaling in $\xx$:
$$
\bm{U}_0(\alpha\xx)=\bm{U}_0(\xx)  \qquad\,\, \mbox{for any $\alpha>0$},
$$
that is, $\bm{U}_0(\xx)$ is constant along the ray originating from $\xx=0$;
in other words, $\bm{U}_0$ depends only on the angular directions of the rays originating from $\xx=0$ in $\mathbb{R}^n$.

A {\it lateral Riemann problem} for \eqref{HCL-1} is a special initial-boundary problem
in a unbounded domain $\mathcal{D}$ that contains the origin and is invariant under the self-similar
scaling ({\it i.e.}, if $\xx\in \mathcal{D}$, then $\alpha \xx\in \mathcal{D}$ for any $\alpha>0$)
so that the initial data and boundary data are also invariant under the self-similar scaling.

Since system \eqref{HCL-1} is invariant under the time-space self-similar scaling, the standard/lateral Riemann
problems are also invariant under the time-space self-similar scaling:
\begin{equation}
(t, \xx)\to (\alpha t, \alpha\xx)\qquad\,\,\mbox{for any $\alpha>0$}.
\end{equation}
Thus, we seek self-similar solutions of the Riemann problems:
\begin{equation}
\bm{U}(t,\xx)=\bm{V}(\frac{\xx}{t}).
\end{equation}
Denote $\xxi=\frac{\xx}{t}$ as the self-similar variables.
Then $\bm{V}(\xxi)$ is determined by
$$
 \D\cdot \bm{F}(\bm{V})-\xxi\cdot \D \bm{V}=0,
$$
that is,
\begin{equation}\label{HCL-1b}
\D\cdot (\bm{F}(\bm{V})-\bm{V}\otimes \xxi)+n\,\bm{V}=0,
\end{equation}
where  $\D=(\partial_{\xi_1}, \cdots, \partial_{\xi_n})$
is the gradient with respect to the self-similar variables
$\xxi=(\xi_1,\cdots, \xi_n)\in \mathbb{R}^n$,
and $\bm{V}\otimes \xxi=(V_i\xi_j)_{1\le i,j\le n}$.
Even though system \eqref{HCL-1} is hyperbolic, system \eqref{HCL-1b} generally
is of mixed elliptic-hyperbolic type, even composite-mixed elliptic-hyperbolic type.
In particular, for a bounded solution $\bm{V}(\xxi)$, system  \eqref{HCL-1b} may be
purely hyperbolic in the far field, {\it i.e.}, outside a large ball in the $\xxi$--coordinates,
but generally is of mixed type or composite-mixed type in a bounded domain containing
the origin, $\xxi=\mathbf{0}$.

For the full Euler system \eqref{HCL-1} with \eqref{HCL-1e}, the self-similar solutions
are governed by the following
system:
\begin{equation}\label{Euler-sms-1}
\left\{\begin{aligned}
&{\rm div}(\rho \vv) + n\rho=0,\\[2mm]
&{\rm div} (\rho \vv\otimes \vv)+ \D p +(n+1)\rho \vv=0,\\
&{\rm div}\big((\frac{1}{2}\rho \lvert \vv\rvert^2+\frac{\gamma p}{\gamma-1})\vv\big)
+n\big(\frac{1}{2}\rho \lvert \vv\rvert^2+\frac{\gamma p}{\gamma-1}\big)=0,
\end{aligned}
\right.
\end{equation}
where $\vv=\uu-\xxi$ is the pseudo-velocity  with
$\bm{V}=(\rho, \rho \vv, \frac{1}{2}\rho \lvert \vv\rvert^2+ \rho e)^\top$.

\smallskip
The weak solutions of system \eqref{HCL-1b} can be defined as follows:

\begin{definition}[Weak Solutions]\label{def-2.1}
A function $\bm{V}\in L_{\rm loc}^\infty(\Lambda)$ in a domain $\Lambda\subset \mathbb{R}^n$ is a weak solution of
system \eqref{HCL-1b} in $\Lambda$, provided that
\begin{equation}
\int_\Lambda \big\{(\bm{F}(\bm{V})-\bm{V}\otimes \xxi)\cdot \D \zeta(\xxi)- n\,\bm{V}\zeta(\xxi)\big\}\, \dd \xxi=0
\quad\,\, \mbox{for any $\zeta\in C_0^1(\Lambda)$}.
\end{equation}
\end{definition}
It can be shown that any weak solution of system \eqref{HCL-1b} in the $\xxi$--coordinates in the sense of Definition \ref{def-2.1}
is a weak solution of system \eqref{HCL-1} in the $(t,\xx)$--coordinates.
Then any co-dimension-one $C^1$--discontinuity $S$
satisfies the Rankine-Hugoniot conditions along $S$ in the $\xxi$--coordinates:
$$
([\bm{F}(\bm{V})]- [\bm{V}]\otimes \xxi)\cdot \nnu_{\rm s}=0,
$$
or equivalently,
\begin{equation}\label{RH-1}
[(\bm{F}(\bm{V})-\bm{V}\otimes \xxi)\cdot \nnu_{\rm s}]=0,
\end{equation}
where $\nnu_{\rm s}$ can be either of the unit normals to $S$, and $[\,\cdot\,]$ denotes the difference between
the traces of the corresponding quantities on the two sides of the co-dimension-one surface $S$.

\smallskip
Similarly, the Riemann
problems for Eq. \eqref{HCL-2} are invariant under the time-space self-similar scaling:
\begin{equation}\label{sss}
(t, \xx,\Phi(t,\xx))\to (\alpha t, \alpha\xx, \frac{\Phi(\alpha t, \alpha\xx)}{\alpha})
\,\,\qquad\mbox{for any $\alpha>0$}.
\end{equation}
Thus, we seek self-similar solutions of the Riemann problem:
\begin{equation}\label{ssform}
\Phi(t,\xx)=t\phi(\frac{\xx}{t}).
\end{equation}
Then $\phi(\xxi)$ is determined by
$$
{\rm div}\,\bm{G}(\phi-\xxi\cdot \D\phi, \D\phi)-\xxi\cdot \D G_0(\phi-\xxi\cdot \D \phi, \D\phi)
=0,
$$
that is,
\begin{equation}\label{HCL-2b}
{\rm div}\big(\bm{G}(\phi-\xxi\cdot \D\phi, \D\phi)- G_0(\phi-\xxi\cdot \D\phi, \D\phi)\xxi\big)
+ n G_0(\phi-\xxi\cdot \D\phi, \D\phi)=0.
\end{equation}
Again, even though Eq. \eqref{HCL-2} is hyperbolic, Eq. \eqref{HCL-2b} generally
is of mixed elliptic-hyperbolic type.
In particular, for a gradient bounded solution $\phi(\xxi)$, Eq. \eqref{HCL-2b} may be
purely hyperbolic in the far field, {\it i.e.}, outside a large ball in the $\xxi$--coordinates,
but generally is of mixed type in a bounded domain containing the origin.

For the Euler equations \eqref{HCL-2} for potential flow with \eqref{1-b1}--\eqref{HCL-2e},
the self-similar solutions are governed by the following second-order quasilinear PDE
for the pseudo-velocity $\varphi=\phi-\frac{1}{2}\lvert \xxi\rvert^2$:
\begin{equation}\label{Euler-PotentialEq-1}
{\rm div}(\rho(\lvert \D\varphi\rvert^2, \varphi)\D\varphi)+ n\rho(\lvert \D\varphi\rvert^2, \varphi)=0,
\end{equation}
where $\rho(\lvert \D\varphi\rvert^2, \varphi)=\big(B_0-(\gamma-1)(\frac{1}{2}\lvert \D\varphi\rvert^2+\varphi)\big)^{\frac{1}{\gamma-1}}$
with $B_0=(\gamma-1)B+1$.

\smallskip
The weak solutions of Eq. \eqref{HCL-2b} can be defined as follows:

\begin{definition}\label{def-2.2}
A function $\phi\in W^{1,\infty}_{\rm loc}(\Lambda)$ in a domain $\Lambda\subset \mathbb{R}^n$ is a weak solution of
system \eqref{HCL-2b} in $\Lambda$, provided that
\begin{align}
\int_\Lambda \big\{
&\big(\bm{G}(\phi-\xxi\cdot \D\phi, \D\phi)- G_0(\phi-\xxi\cdot \D\phi, \D\phi)\xxi\big)\cdot \D\zeta(\xxi)
\nonumber\\
&-n G_0(\phi-\xxi\cdot \D\phi, \D\phi)\zeta(\xxi)\big\}\, \dd \xxi=0
 \label{2.16a}
\end{align}
for any $\zeta\in C_0^1(\Lambda)$.
\end{definition}
Similarly, it can shown that any weak solution of Eq. \eqref{HCL-2b} in the $\xxi$--coordinates in the sense of Definition \ref{def-2.2}
is a weak solution of Eq. \eqref{HCL-2} in the $(t,\xx)$--coordinates. Then any co-dimension-one
$C^1$--discontinuity $S$
satisfies the Rankine-Hugoniot conditions along $S$ in the $\xxi$--coordinates:
$$
[\phi]=0, \qquad [\bm{G}(\phi-\xxi\cdot \D\phi, \D\phi)- G_0(\phi-\xxi\cdot \D\phi, \D\phi)\xxi]\cdot\nnu_{\rm s}=0,
$$
or equivalently,
$$
[\phi]=0, \qquad [(\bm{G}(\phi-\xxi\cdot \D\phi, \D\phi)- G_0(\phi-\xxi\cdot \D\phi, \D\phi)\xxi)\cdot\nnu_{\rm s}]=0,
$$
where $\nnu_{\rm s}$ is either of the unit normals to $S$.

\section{Two-Dimensional Riemann Problem I: \\
Two Shocks and Two Vortex Sheets \\ for the Pressure Gradient System}

In this section, we present the first 2-D Riemann problem, Riemann Problem I,
through the pressure gradient system that is
a hyperbolic system of conservation laws.

The pressure gradient system takes the following form:
\begin{equation}\label{pgs}
\begin{cases}
u_t+p_{x_1}=0,\\[2pt]
v_t+p_{x_2}=0,\\[2pt]
E_t+(pu)_{x_1}+(pv)_{x_2}=0,
\end{cases}
\end{equation}
where  $E=\frac{\lvert \uu\rvert^2}{2}+p$ with $\uu=(u,v)$.
System \eqref{pgs} can be written in form \eqref{HCL-2} with
\begin{equation}\label{pgs-b}
\bm{U}=(\uu,E)^\top,
\,\,\, \FF_1=(E-\frac{{\lvert \uu\rvert^2}}{2})(1, 0, u)^\top,
\,\,\, \FF_2=(E-\frac{\lvert \uu\rvert^2}{2})(0, 1, v)^\top.
\end{equation}

There are two mechanisms for the fluid motion: the inertia and the pressure differences.
Corresponding to a separation of these two mechanisms,
the full Euler equations \eqref{HCL-1} with \eqref{HCL-1e} in gas dynamics
can be split into two subsystems of conservation laws:
the pressure gradient system and the pressureless Euler system, respectively;
also see \cite{agarwal1994modified,ChenLeFloch,li1985second} and the references cited therein
for this and
similar flux-splitting ideas which have been widely used in order to design the so-called flux-splitting
schemes and their high-order accurate extensions.
Furthermore, system \eqref{pgs} can also be deduced from
system  \eqref{HCL-1} with \eqref{HCL-1e} under the physical regime whereby
the velocity
is small and the adiabatic gas constant $\gamma$ is large; see Zheng \cite{yuxi1997existence}.
An asymptotic derivation of system \eqref{pgs} has also been presented by Hunter as described in \cite{zheng2006two}.
We refer the reader to \cite{li1998two,zheng2012systems} for further background on system \eqref{pgs}.

\subsection{2-D Riemann Problem I:\\ Two Shocks and Two Vortex Sheets}

\begin{figure}[htbp]
\vspace{-1.1cm}
\hspace{-1.35cm}
\includegraphics[scale=1.2]{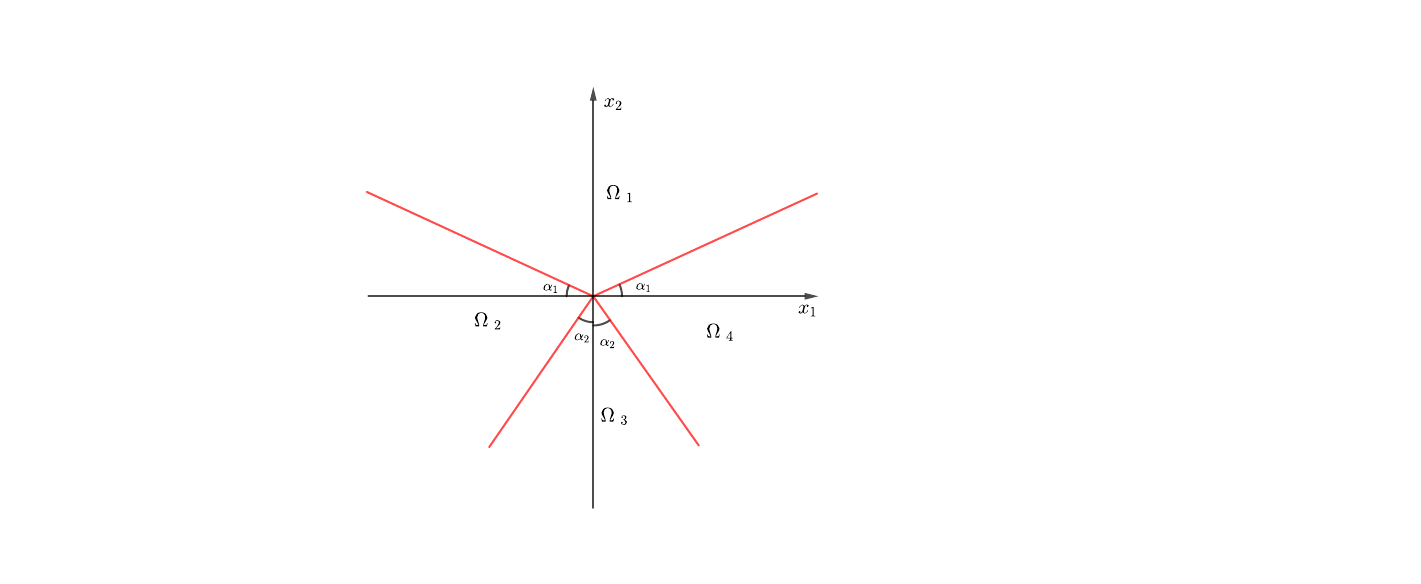}
\setlength{\abovecaptionskip}{-1.0cm}
\caption{Riemann Problem I: Riemann initial data ({\it cf}. \cite{CWZ,zheng2003global})}\label{figure1}
\end{figure}

We now consider the following Riemann problem:

\begin{problem}[2-D Riemann Problem I: Two Shocks and Two Vortex Sheets]\label{prob:3.1}
Seek a global solution of system \eqref{pgs} with Riemann initial data that consist of four constant states
in four sectorial regions
$\Omega_i$  with symmetric sectorial angles $($see {\rm Fig. \ref{figure1}}$)${\rm :}
\begin{equation}\label{initialdata}
(p,\uu)(0,\xx)=(p_i,\uu_i) \qquad\,\, \mbox{for $\,\,\xx\in \Omega_i, \, i=1,2,3,4$},
\end{equation}
such that the four initial constant states are required to satisfy the following conditions{\rm :}
\begin{equation}\label{initialdate1}
\begin{cases}
\text{A forward shock $S_{41}^{+}$ is formed between states $(1)$ and $(4)$,}\\[2pt]
\text{A backward shock $S_{12}^{-}$ is formed between states $(1)$ and $(2)$,}\\[2pt]
\text{A vortex sheet $J_{23}^{+}$ is formed between states $(2)$ and $(3)$,}\\[2pt]
\text{A vortex sheet $J_{34}^{-}$ is formed between states $(3)$ and $(4)$}.
\end{cases}
\end{equation}
\end{problem}

\vspace{-10pt}

This Riemann problem initially with the assumption that angle $\alpha_1=\alpha_2$ is close to zero
was first analyzed rigorously in Zheng \cite{zheng2003global},
for which the two shocks bend slightly and the diffracted shock $\Gamma_{\rm shock}$ does
not meet the inner sonic circle $C_2$.
In the recent work \cite{CWZ}, this Riemann problem has been solved globally for the general case; that is,
the angle between the two shocks is not necessarily close to $\pi$.

\subsection{Reformulation of Riemann Problem I}

As discussed earlier,
we seek self-similar solutions in the self-similar coordinates with the form:
\begin{equation*}
(p,\uu)(t,\xx)=(p,\uu)(\xxi) \qquad \text{with $\xxi=\frac{\xx}{t},\,t>0$}.
\end{equation*}
In the $\xxi$--coordinates, system \eqref{pgs} can be rewritten in form
\eqref{HCL-1b} with \eqref{pgs-b}.
The four waves in Riemann Problem I can be obtained by solving four 1-D Riemann problems in the self-similar
coordinates $\xxi$, which form the following configuration
as shown in Fig. \ref{figure2}:

\begin{figure}[htbp]
\vspace{-0.3cm}
\includegraphics[scale=0.9]{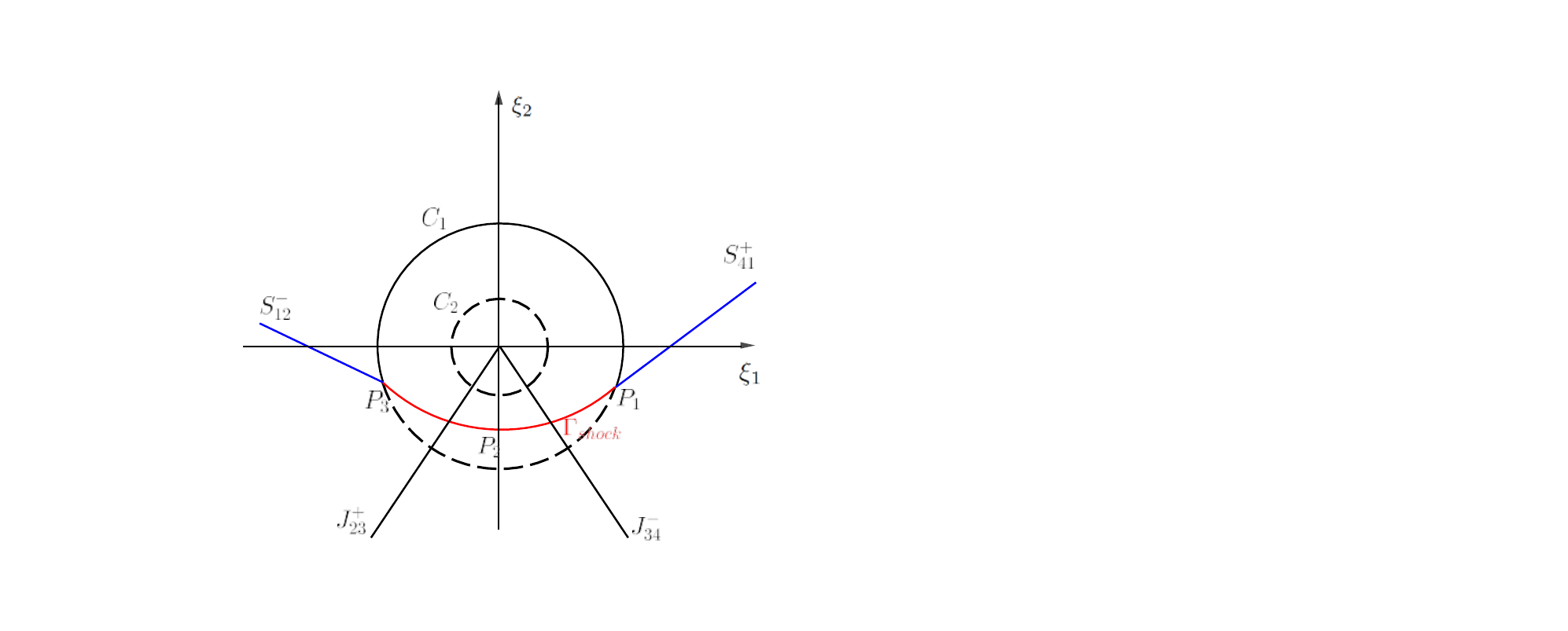}
\setlength{\abovecaptionskip}{-1.2cm}
\caption{Riemann Problem I:  Riemann solution configuration ({\it cf.} \cite{CWZ})}\label{figure2}
\end{figure}

More precisely, let $\xi_2=f(\xi_1)$ be a $C^1$--discontinuity curve of a bounded discontinuous solution of
system \eqref{HCL-1b} with \eqref{pgs-b}.
From the Rankine-Hugoniot relations on $\xi_2=f(\xi_1)$:
\begin{equation*}
\begin{cases}
\big(\xi_1f'(\xi_1)-f(\xi_1)\big)[u]-f'(\xi_1)[p]=0,\\[1mm]
\big(\xi_1 f'(\xi_1)-f(\xi_1)\big)[v]+[p]=0,\\[1mm]
\big(\xi_1 f'(\xi_1)-f(\xi_1)\big)[E]-f'(\xi_1)[pu]+[pv]=0,
\end{cases}
\end{equation*}
we find that $\xi_2=f(\xi_1)$ can be one of the two nonlinear discontinuities:
\begin{equation}\label{shock11}
\begin{cases}
\frac{\dd f(\xi_1)}{\dd\xi_1}=\sigma_{\pm}=-\frac{[u]}{[v]}
  =\frac{\xi_1 f(\xi_1)\pm\sqrt{\overline{p}( \xi_1^2+  \lvert f(\xi_1)\rvert^2 -\overline{p})}}{\xi_1^2-\overline{p}},\\[1.5mm]
[p]^2=\overline{p}([u]^2+[v]^2),
\end{cases}
\end{equation}
or a vortex sheet (linearly degenerate discontinuity):
\begin{equation}\label{shock12}
\begin{cases}
\sigma_{0}=\frac{f(\xi_1)}{\xi_1}=\frac{[v]}{[u]},\\[5pt]
[p]=0,
\end{cases}
\end{equation}
where $\overline{p}$
is the average of the pressure
on the two sides of the discontinuity.

A nonlinear discontinuity is called a shock if it satisfies \eqref{shock11} and the entropy
condition: {\it pressure $p$ increases across it in the flow direction}; that is,
the pressure ahead of the wave-front is larger than that behind the wave-front.
There are two types of shocks $S^{\pm}$:

\smallskip
\begin{itemize}
\item$ S=S^{+}$ if $\D p$ and the flow direction form a right-hand system;

\smallskip
\item $S=S^{-}$ if $\D p$ and the flow direction form a left-hand system.
\end{itemize}
A discontinuity is called a vortex sheet if it satisfies \eqref{shock12}.
There are two types of vortex sheets $J^{\pm}$ determined by the signs of the vorticity:
$$
J^{\pm}: \quad \text{curl}\,\uu=\pm\infty.
$$

It can be shown that,
for fixed $(p_1,\uu_1)$ and $p_2=p_3=p_4$ satisfying $p_1>p_2$,
there exist states $\uu_i, i=2,3,4$, depending on
angles $(\alpha_1,\alpha_2)$ continuously such that the conditions
in \eqref{initialdate1} for the Riemann initial data hold.

\begin{figure}[htbp]
\vspace{-0.60cm}
\hspace{-2.2cm}
{\includegraphics[scale=0.80]{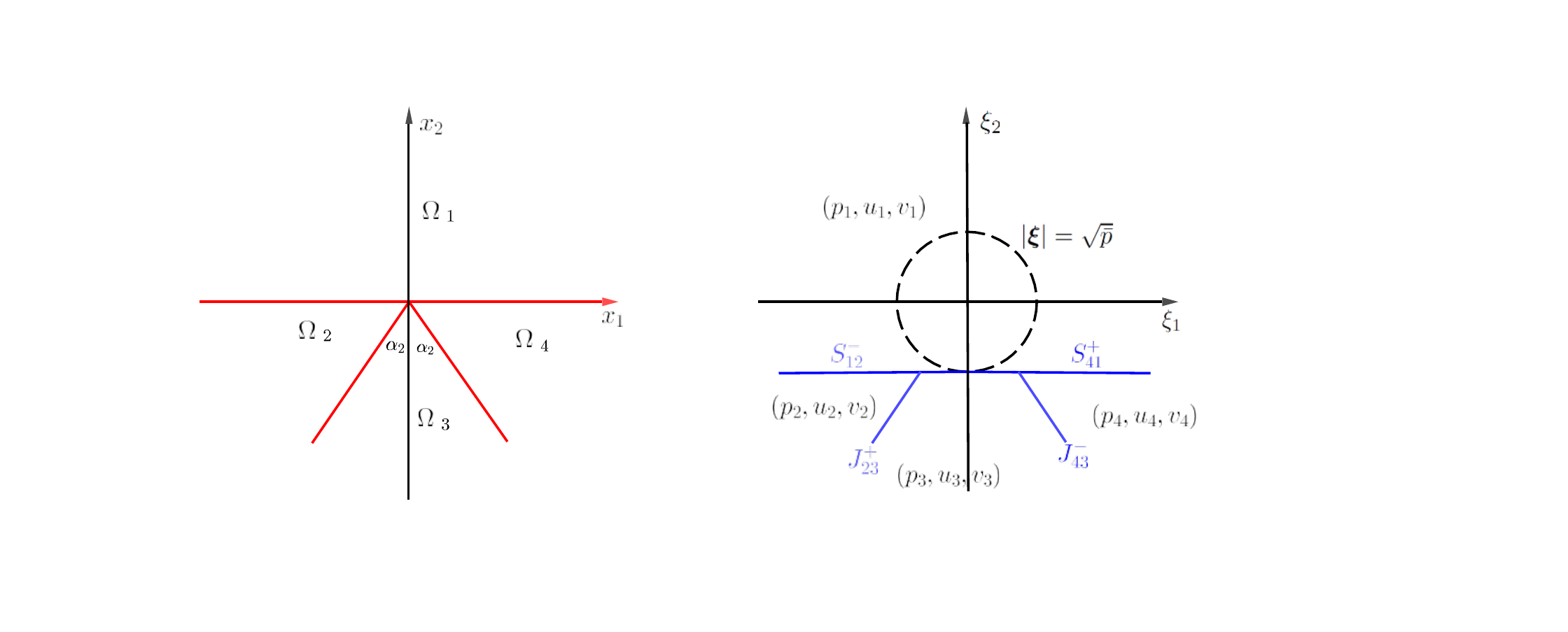}\qquad\qquad}
{\setlength{\abovecaptionskip}{-1.3cm}\caption{\label{figure3}
{The Riemann data and the global solution when $\alpha_1=0$ ({\it cf}. \cite{CWZ})}}}
\end{figure}

There is a critical case when $\alpha_1=0$.
Then the Riemann initial data satisfy
$$
p_1>p_2=p_3=p_4, \quad u_1=u_2=u_3=u_4,\quad v_1>v_2=v_3=v_4.
$$
The global Riemann solution is a piecewise constant solution with two planar shocks:
$S_{12}^{-}$ for $\xi_1<0$ and $S_{41}^{+}$ for $\xi_1>0$ on the line:
$\xi_2=-\sqrt{\overline{p}}$, with
$$
[v]=-\frac{[p]}{\sqrt{\overline{p}}}, \quad\, [u]=0
\qquad\,\,\, \mbox{for  $\overline{p}=\frac{p_1+p_2}{2}$},
$$
and two vortex sheets $J_{23}^{+}$ and $J_{34}^{-}$,
as shown in Fig. \ref{figure3}.
The two planar shocks $S_{12}^{-}$ and $S_{41}^{+}$ are both tangential to the circle,
$\lvert \xxi\rvert=\sqrt{\overline{p}}$,
with the tangent point on the circle as the end-point.
It follows from the expression of $J_{23}^{+}$ given in \eqref{shock12}
that $p_2=p_3$ on both sides of $J_{23}^{+}$.
At the point where $J_{23}^{+}$ intersects with $S_{12}^{-}$,
we see that $J_{23}^{+}$ does not affect the shock owing to $p_2=p_3$.
The intersection between $J_{34}^{-}$ and $S_{41}^{+}$ can be handled
in the same way.

\medskip
We now consider the general case: $\alpha_1\in (0,\frac{\pi}{2})$.
From system \eqref{HCL-1b} with \eqref{pgs-b},
we can derive the following second-order nonlinear equation for $p$:
\begin{equation}\label{eqforp}
(p-\xi_1^{2})p_{\xi_1\xi_1}-2\xi_1\xi_2 p_{\xi_1\xi_2}+(p-\xi_2^{2})p_{\xi_2\xi_2}
+\frac{(\xi_1 p_{\xi_1}+\xi_2 p_{\xi_2})^{2}}{p}
-2(\xi_1 p_{\xi_1}+\xi_2 p_{\xi_2})=0.
\end{equation}
Eq. \eqref{eqforp} is of mixed hyperbolic-elliptic
type, which is hyperbolic when $\lvert\xxi\rvert>\sqrt{p}$ and elliptic when $\lvert\xxi\rvert<\sqrt{p}$
with the transition boundary -- the sonic circle $\lvert\xxi\rvert=\sqrt{p}$.
Furthermore, in the polar coordinates:
$
(r,\theta)=(\lvert \xxi\rvert, \arctan(\frac{\xi_2}{\xi_1})),
$
Eq. \eqref{eqforp} becomes
\begin{equation}
\label{eqforp2}
Qp:=(p-r^{2})p_{rr}+\frac{p}{r^{2}}p_{\theta\theta}+\frac{p}{r}p_{r}+\frac{1}{p}\big((r p_{r})^{2}-2r p p_{r}\big)=0,
\end{equation}
which is  hyperbolic when $p<r^2$ and elliptic when $p>r^2$.
The sonic circle is given by $r=r(\theta)=\sqrt{p(r(\theta),\theta)}$.

In the $\xxi$--coordinates,
the four waves come from the far-field (at infinity, corresponding to $t=0$)
and keep planar waves before the two shocks meet the outer sonic circle $C_1$
of state $(1)$:
$$
C_1:=\{\,\xxi\,:\,\lvert\xxi\rvert=\sqrt{p_1}\,\}.
$$
When the two shocks $S_{12}^{-}$ and $S_{41}^{+}$ meet the sonic circle $C_1$
at points $P_3$ and $P_1$ respectively,
the key issue is whether they bend and meet to
form a diffracted shock, denoted by $\Gamma_{\rm shock}$;
see Fig. \ref{figure2}.
Since the whole configuration is symmetric with respect to the
$\xi_2$--axis,
$\Gamma_{\rm shock}$ must be perpendicular to $\xi_1=0$ at point $P_2$
where the two diffracted shocks meet.
It is not known {\it a priori} whether the diffracted shock may
degenerate partially into a portion of the inner sonic circle $C_2$ of state $(2)$.
Once this case occurs, $p=p_2$ on the sonic circle, which satisfies the oblique derivative
condition on the diffracted shock automatically.
Observe that the two vortex sheets
$J_{23}^{+}$ and $J_{34}^{-}$ and the diffracted shock $\Gamma_{\rm shock}$
have no influence on each other during the intersection,
as pointed out earlier by Zhang-Li-Zhang \cite{ZLZ}.
Therefore, from now on, we first {\it ignore} the two vortex sheets
and focus mainly on the diffracted shock.

\smallskip
On $\Gamma_{\rm shock}$, the Rankine-Hugoniot conditions in the polar coordinates
must be satisfied:
\begin{equation}\label{eq:2.1}
\begin{cases}
r[u]-\big(\cos\theta +\frac{1}{r}\frac{\dd r}{\dd\theta}\sin\theta\big)[p]=0,\\[1mm]
r[v]-\big(\sin\theta -\frac{1}{r}\frac{\dd r}{\dd\theta}\cos\theta\big)[p]=0,\\[1mm]
r[E]-\big(\cos\theta+\frac{1}{r}\frac{\dd r}{\dd\theta}\sin\theta\big)[pu]
  -\big(\sin\theta - \frac{1}{r}\frac{\dd r}{\dd\theta}\cos\theta\big)[pv]=0.
\end{cases}
\end{equation}
Owing to
$
[pu]=\overline{p}\,[u]+\overline{u}\,[p],
$
with $\overline{p}$ as the average of the two neighboring
states of $p$,
we eliminate $[u]$ and $[v]$ in the third equation in \eqref{eq:2.1} to obtain
\begin{equation*}
\Big(\frac{\dd r}{\dd\theta}\Big)^2=\frac{r^2(r^2-\overline{p})}{\overline{p}}.
\end{equation*}

The shock diffraction can also be regarded to be generated from point $P_2$ in two directions,
which implies that $r^{\prime}(\theta)>0$ for $\theta\in [\frac{3\pi}{2},\theta_{1}]$
and $r^{\prime}(\theta)<0$ for $\theta\in [\theta_{3},\frac{3\pi}{2}]$,
where $\theta_{i}$ are denoted as the $\theta$--coordinates of points $P_i$, $i=1,3$,
respectively.
Thus, we choose
\begin{equation}\label{shockeq}
\frac{\dd r}{\dd\theta}=g(p(r(\theta),\theta),r(\theta)):=
\begin{cases}
r\sqrt{\frac{r^2-\overline{p}}{\overline{p}}}\qquad\,\, &\text{for}\,\, \theta\in [\frac{3\pi}{2},\theta_{1}],\\[5pt]
-r\sqrt{\frac{r^2-\overline{p}}{\overline{p}}}\qquad\,\,  &\text{for}\,\, \theta\in [\theta_{3},\frac{3\pi}{2}].
\end{cases}
\end{equation}
It follows from \eqref{shock11}, or \eqref{eq:2.1}, that
$
[p]^2=\overline{p}\,\big([u]^2+[v]^2\big).
$
Then taking the derivative $r^\prime(\theta)\partial_r+\partial_\theta$ on both sides of this equation
along the shock  yields the derivative boundary condition
on $\Gamma_{\rm shock}=\{(r(\theta),\theta)\,:\, \theta_3\le \theta\le \theta_1\}$:
\begin{equation}\label{eqrh}
\beta_1 p_r+\beta_2 p_\theta=0,
\end{equation}
where $\mathbf{\beta}=(\beta_1,\beta_2)$ is a function of $(p,p_2,r(\theta),r^\prime(\theta))$ with
\begin{equation}
\label{eqrh2}
\beta_1=2r^\prime(\theta)\Big(\frac{r^2-\overline{p}}{r^2}-\frac{[p]}{4\overline{p}}
+\frac{\overline{p}(r^2-p)}{r^2p}\Big), \quad
\beta_2=\frac{4(r^2-\overline{p})}{r^2}-\frac{[p]}{2\overline{p}}.
\end{equation}
The obliqueness becomes
$$
\mu:=(\beta_1,\beta_2)\cdot (1,-r^\prime(\theta))=-2r^\prime(\theta)\big(1-\frac{\overline{p}}{p}\big).
$$
Note that $\mu$ vanishes at point $P_2$  where $r^\prime(\frac{3\pi}{2})=0$ and
$$
\beta_1=0,\qquad \beta_2=-\frac{[p]}{2\overline{p}}<0,
$$
owing to  $p>p_2$.

Let $\Gamma_{\rm sonic}$ be the larger portion $\widehat{P_1P_3}$ of the sonic circle $C_1$
of state $(1)$. On $\Gamma_{\rm sonic}$, $p$ satisfies the Dirichlet boundary condition:
\begin{equation}\label{dirichlet}
p=p_1.
\end{equation}
Let $\Omega$ be the bounded domain enclosed by $\Gamma_{\rm sonic}$
and $\Gamma_{\rm shock}$.
Then Riemann Problem I (Problem 3.1)
can be reformulated into
the following free boundary problem:

\begin{problem}[Free Boundary Problem]\label{problem}
Seek a solution $(p(r,\theta), r(\theta))$
such that $p(r, \theta)$ and $r(\theta)$ are determined by
Eq. \eqref{eqforp2} in $\Omega$ and
the free boundary conditions \eqref{shockeq}--\eqref{eqrh2} on $\Gamma_{\rm shock}$
$($the derivative boundary condition$)$,
in addition to
the Dirichlet boundary condition \eqref{dirichlet} on $\Gamma_{\rm sonic}$.
\end{problem}

\subsection{Global Solutions of Riemann Problem I: \\Free Boundary Problem, Problem \ref{problem}}\label{sec2.3}

To solve Riemann Problem I, it suffices to deal with the free boundary problem, Problem \ref{problem},
which has been solved as stated in the following theorem.

\begin{theorem}[Chen-Wang-Zhu \cite{CWZ}]\label{thm1}
There exists a global solution $(p(r,\theta), r(\theta))$ of {\rm Problem \ref{problem}}
in domain $\Omega$ with the free boundary
$$
\Gamma_{\rm shock}:=\{(r(\theta),\theta)\,:\, \theta_{3}\le \theta\le \theta_{1}\}
$$
such that
\begin{equation*}
p\in C^{2,\alpha}(\Omega)\cap C^{\alpha}(\overline{\Omega}), \qquad
r\in C^{2,\alpha}((\theta_{3},\theta_{1}))\cap C^{1,1}([\theta_{3},\theta_{1}]),
\end{equation*}
where $\alpha\in (0,1)$ depends only on the Riemann initial data.
Moreover, the global solution $(p(r,\theta),r(\theta))$
satisfies the following properties{\rm :}

\smallskip
\begin{itemize}
\item[\emph{(i)}] $p>p_2$ on the free boundary $\Gamma_{\rm shock}${\rm ;} that is,
$\Gamma_{\rm shock}$ does not meet the sonic circle $C_2$ of state $(2)$.

\item[\emph{(ii)}] The free boundary $\Gamma_{\rm shock}$ is
convex in the self-similar coordinates.

\item[\emph{(iii)}] The global solution $p(r,\theta)$ is $C^{1,\alpha}$ up to the sonic boundary
$\Gamma_{\rm sonic}$ and Lipschitz continuous across $\Gamma_{\rm sonic}$.

\item[\emph{(iv)}] The Lipschitz regularity of the solution across $\Gamma_{\rm sonic}$
from the inside of the subsonic domain is optimal.
\end{itemize}
\end{theorem}

There are three main difficulties for the proof of Theorem \ref{thm1}:
\smallskip
\begin{enumerate}[(i)]
\item The diffracted shock $\Gamma_{\rm shock}$ is a free boundary, which is not known {\it a priori} whether
it coincides with the inner sonic circle $C_2$ of state $(2)$.
\item  On the sonic boundary $\Gamma_{\rm sonic}$, owing to $p_1=r^2$, the ellipticity of Eq. \eqref{eqforp2} degenerates.
\item At point $P_2$ where the diffracted shock $\Gamma_{\rm shock}$ meets the $\xi_2$--axis: $\xi_1=0$,
    the obliqueness of derivative boundary condition fails, since
     $$
      (\beta_1,\beta_2)\cdot (1,-r^{\prime}(\theta))=0.
     $$

\end{enumerate}

\smallskip
In the proof of  Theorem \ref{thm1},
we first assume that $p\geq p_2+\delta$ holds on $\Gamma_{\rm shock}$ for some $\delta>0$;
that is, $\Gamma_{\rm shock}$ cannot coincide with the sonic circle $C_2$
of state $(2)$, which is eventually proved.
For the third difficulty,  we may express this as a one-point Dirichlet condition $p(P_2)=\hat{p}$
by solving
$$
2r(\theta_2)=p(r(\theta_2),\theta_2)+p_2.
$$
More precisely, the existence proof is divided into four steps:

\medskip
1. Since Eq. \eqref{eqforp2} degenerates on the sonic boundary, the differential operator $Q$
in Eq. \eqref{eqforp2}
is replaced by the regularized operator:
$$
Q^{\varepsilon}=Q+\varepsilon\Delta_\xxi.
$$
The free boundary $\Gamma_{\rm shock}$ is first fixed,
then the equation and the derivative boundary condition are linearized,
and the existence of a solution of the linear fixed mixed-type boundary problem
is established for the regularized equation in the polar coordinates.

\smallskip
2. Based on the estimates of solutions of the linear fixed boundary problem obtained in Step $1$,
the existence of a solution of the nonlinear fixed boundary problem is proved via
the Schauder fixed point theorem.

\smallskip
3. The existence of a solution of the free boundary problem
with the oblique derivative boundary condition for the regularized elliptic equation
is established by using the Schauder fixed point argument again.
It follows that the free boundary  never meets the sonic circle $C_2$
of state $p_2$.

\smallskip
4. Finally, the limiting solution as the elliptic regularization parameter $\varepsilon$ tends to $0$
is proved to be a solution of Problem \ref{problem}.

\medskip
In Theorem \ref{thm1}, a global solution $p$ of the second-order
equation \eqref{eqforp} in $\Omega$ is constructed, which is piecewise
constant in the supersonic domain.
Moreover, $p$ is proved to be Lipschitz continuous across the degenerate
sonic boundary $\Gamma_{\rm sonic}$ from $\Omega$ to the supersonic domain.
To recover velocity $\uu=(u, v)$, we consider the first two equations in
system \eqref{HCL-1b} with \eqref{pgs-b}.
We can rewrite these equations in the radial variable $r$ as
$$
\frac{\partial \uu}{\partial r}=\frac{1}{r}\D p,
$$
and integrate from the boundary of the subsonic domain toward the origin.
It is direct to see that $\uu$ is at least Lipschitz continuous across $\Gamma_{\rm sonic}$.
Furthermore, $\uu$ has the same regularity as $p$ inside $\Omega$ except origin $r=0$.
However, $\uu$ may be multi-valued at the origin ({\it i.e.}, $r=0$).  Therefore, we have

\begin{theorem}[Chen-Wang-Zhu \cite{CWZ}]\label{thm2}
Let the Riemann initial data satisfy \eqref{initialdate1}.
Then there exists a global solution $(p, \uu)(r,\theta)$ with the $2$-D shock
$$
\Gamma_{\rm shock}=\{(r(\theta), \theta)\,:\, \theta_{3}\le\theta\le \theta_{1}\}
$$
such that
\begin{equation*}
(p, \uu)\in C^{2,\alpha}(\Omega), \quad p\in C^{\alpha}(\overline{\Omega}), \quad
r\in C^{2,\alpha}((\theta_{3},\theta_{1}))\cap C^{1,1}([\theta_{3},\theta_{1}]),
\end{equation*}
and $(p, \uu)$ are piecewise constant in the supersonic domain.
Moreover, the global solution $(p, \uu)$ with shock $\Gamma_{\rm shock}$
satisfies properties {\rm (i)}--{\rm (ii)} in {\rm Theorem \ref{thm1}} and

\smallskip
\emph{(a)} $(p, \uu)$ is $C^{\alpha}$ up to the sonic boundary $\Gamma_{\rm sonic}$ and Lipschitz continuous across $\Gamma_{\rm sonic}$.

\emph{(b)} The Lipschitz regularity of both solution $(p, \uu)$ across $\Gamma_{\rm sonic}$ from the subsonic domain $\Omega$
and shock $\Gamma_{\rm shock}$ across points $\{P_1, P_3\}$
is optimal.
\end{theorem}

\smallskip
More details can be found in Chen-Wang-Zhu \cite{CWZ}. Similar results can be obtained for the nonlinear wave system introduced
in Section 4 below by using the same approach and related techniques/methods.
Furthermore, Riemann Problem I for the Euler equations for potential flow has also been solved recently in \cite{CCLW}.

\section{Two-Dimensional Riemann Problem II: \\
The Lighthill Problem for Shock Diffraction for the Nonlinear Wave System}

In this section, we present the second Riemann problem,
Riemann problem II -- the Lighthill problem for shock diffraction by 2-D convex cornered wedges
in compressible fluid flow (Lighthill \cite{Lighthill1,Lighthill2}),
through the nonlinear wave system; also see \cite{Bargman,CDX,FWB,FTB}.

The nonlinear wave system consists of three conservation laws, which
takes the form:
\begin{equation} \label{1.1a}
\begin{cases}
\rho_{t}+m_{x_1}+n_{x_2}=0,\\
m_{t}+p_{x_1}=0,\\
n_{t}+p_{x_2}=0,
\end{cases}
\end{equation}
for $(t,\mathbf{x})\in[0,\infty)\times \mathbb{R}^2$,
where $\rho$ stands for the density,  $p$ for the
pressure, and $(m,n)$ for the momenta in the $\xx$--coordinates.
The pressure-density constitutive relation is
\begin{equation}\label{pressure}
p(\rho)=\frac{\rho^\gamma}{\gamma} \qquad\,\,\mbox{for $\gamma>1$}
\end{equation}
by scaling without loss of generality. Then the sonic speed
$c=c(\rho)$ is determined by
$$
c(\rho):=\sqrt{p'(\rho)}=\rho^{\frac{\gamma-1}{2}},
$$
which is a positive, increasing function for all $\rho>0$.
System \eqref{1.1a} can be written in form \eqref{HCL-1} with
\begin{equation}\label{1.1ab}
\bm{U}=(\rho, m, n)^\top, \quad
\bm{F}_1=(m, p,0)^\top, \quad
\bm{F}_2=(n, 0, p)^\top.
\end{equation}

The 2-D nonlinear wave system \eqref{1.1a} is derived
from the compressible isentropic gas dynamics by neglecting the
inertial terms, {\it i.e.}, the quadratic terms in the velocity; see
Canic-Keyfitz-Kim \cite{canic2006free}.

\subsection{Riemann Problem II: The Lighthill Problem for Shock Diffraction by Convex Cornered Wedges}

Let $S_0$ be the vertical planar shock in the
$(t, {\bf x})$--coordinates,
with the left constant state
$U_{1}=(\rho_{1},m_{1},0)$ and the right state
$U_{0}=(\rho_{0},0,0)$, satisfying
$$
m_1=\sqrt{\big(p(\rho_{1})-p(\rho_{0})\big)(\rho_{1}-\rho_{0})}>0,
\qquad \rho_1>\rho_0.
$$
When $S_0$ passes through a convex cornered wedge:
$$
W:=\{\xx=(x_1, x_2)\, :\, x_2<0,  x_1\le x_2\,{\rm
ctan}\,\theta_{\rm w}\},
$$
shock diffraction occurs, where the wedge angle $\theta_{\rm w}$ is
between $-\pi$ and $0$; see Fig. \ref{cog}. Then the shock
diffraction problem can be formulated as follows:

\begin{figure}[h]
 \vspace{-0.2cm}
  \centering
  \includegraphics[width=0.47\textwidth]{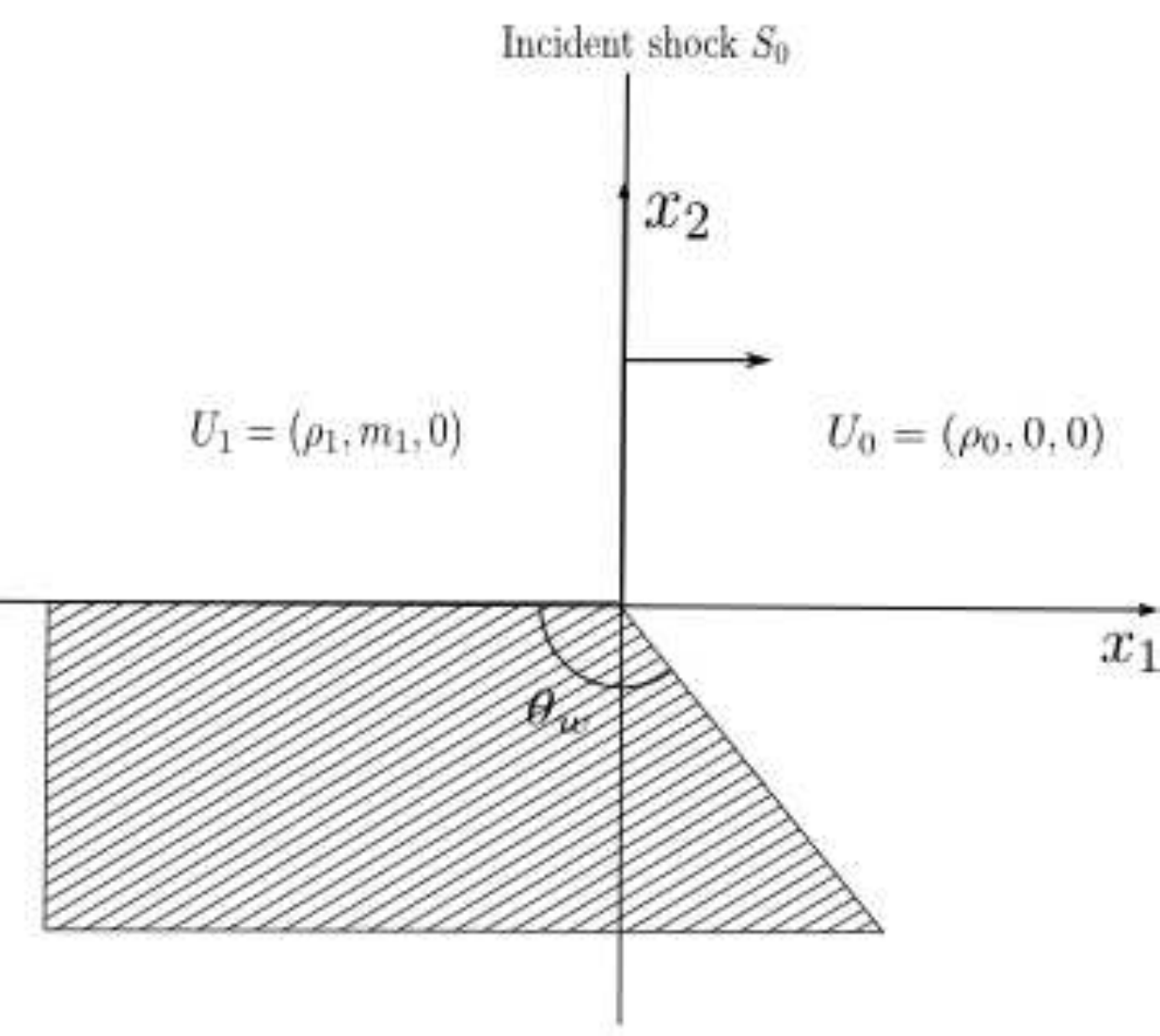} \quad
  \includegraphics[width=0.46\textwidth]{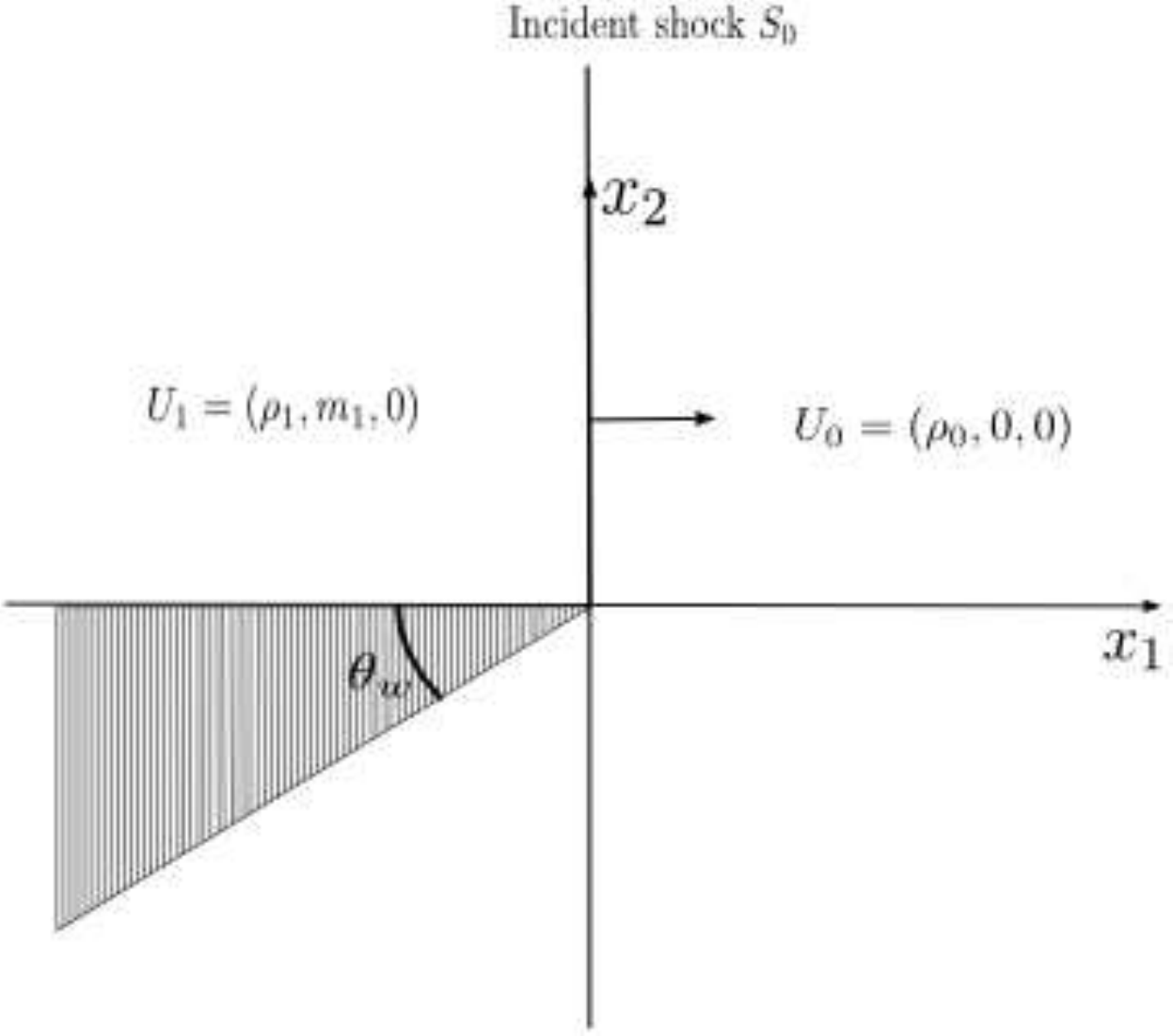}
  \caption{Riemann Problem II: The Lighthill problem ({\it cf.} \cite{CDX})}.
  \label{cog}
\end{figure}

\vspace{-18pt}
\begin{problem}[Riemann Problem II: The Lighthill Problem for Shock Diffraction]\label{prob:4.1}
 {\it Seek a
solution of system \eqref{1.1a}--\eqref{pressure} with the initial condition at $t=0${\rm :}
\begin{equation}\label{initial-condition}
\bm{U}\rvert_{t=0}
=\begin{cases}
(\rho_0, 0, 0) \quad &\mbox{\rm in $\{-\pi+\theta_{\rm w}\le \arctan(\frac{x_2}{x_1})\le \frac{\pi}{2}\}$},\\
(\rho_1, m_1, 0) \quad &\mbox{\rm in $\{x_1<0, x_2>0\}$},
\end{cases}
\end{equation}
and the slip boundary condition along the wedge boundary $\partial
W${\rm :}
\begin{equation}\label{boundary-condition}
(m, n)\cdot \nnu_{\rm w}\mid_{\partial W}=0,
\end{equation}
where $\nnu_{\rm w}$ is the exterior unit normal to $\partial W$ $($see Fig.
{\rm \ref{cog}}$)$.}
\end{problem}

\subsection{Reformulation of Riemann Problem II}

\begin{figure}[h]
\vspace{-5pt}
  \centering
  \includegraphics[width=0.47\textwidth]{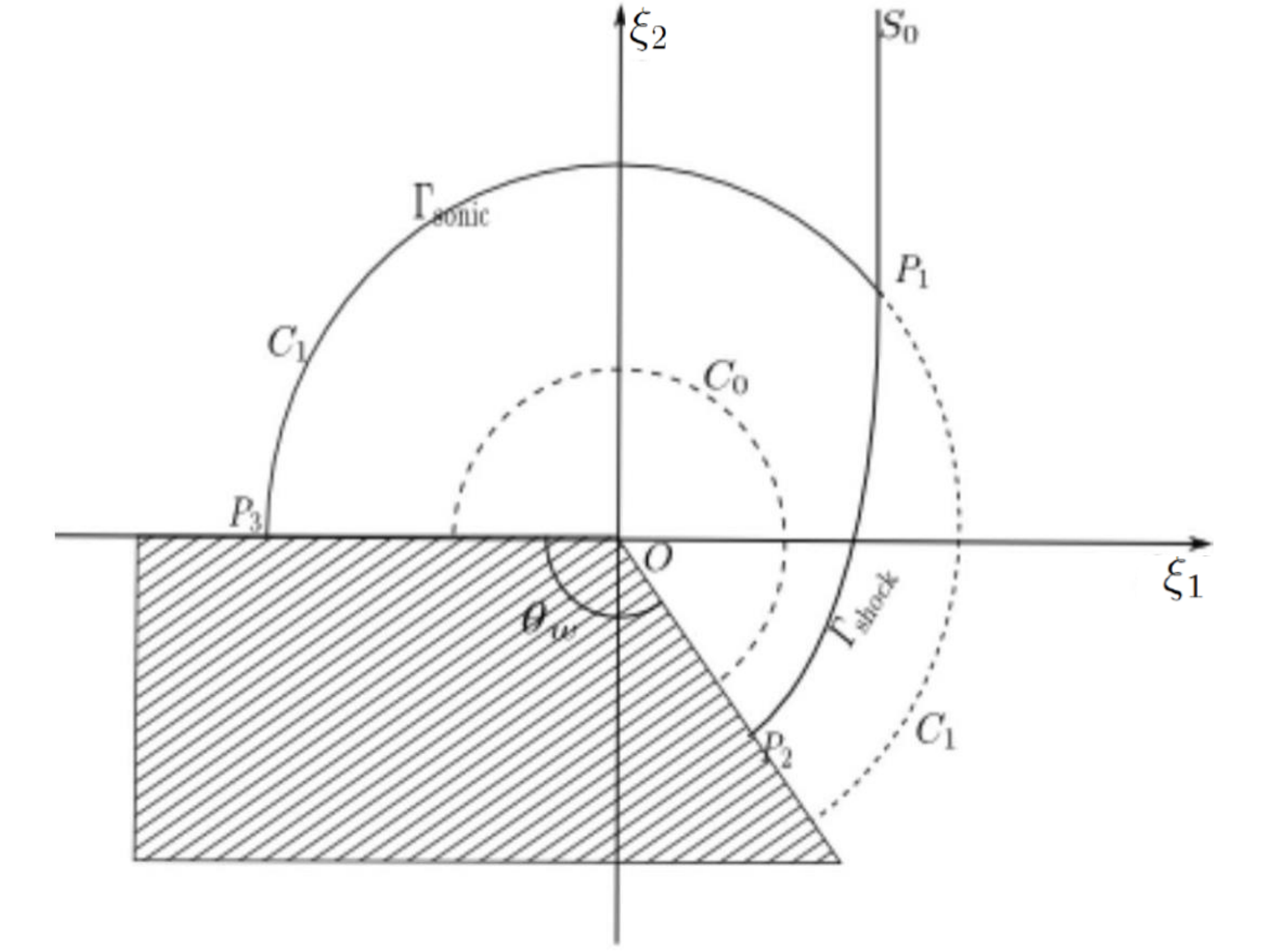}\,\,
  \includegraphics[width=0.47\textwidth]{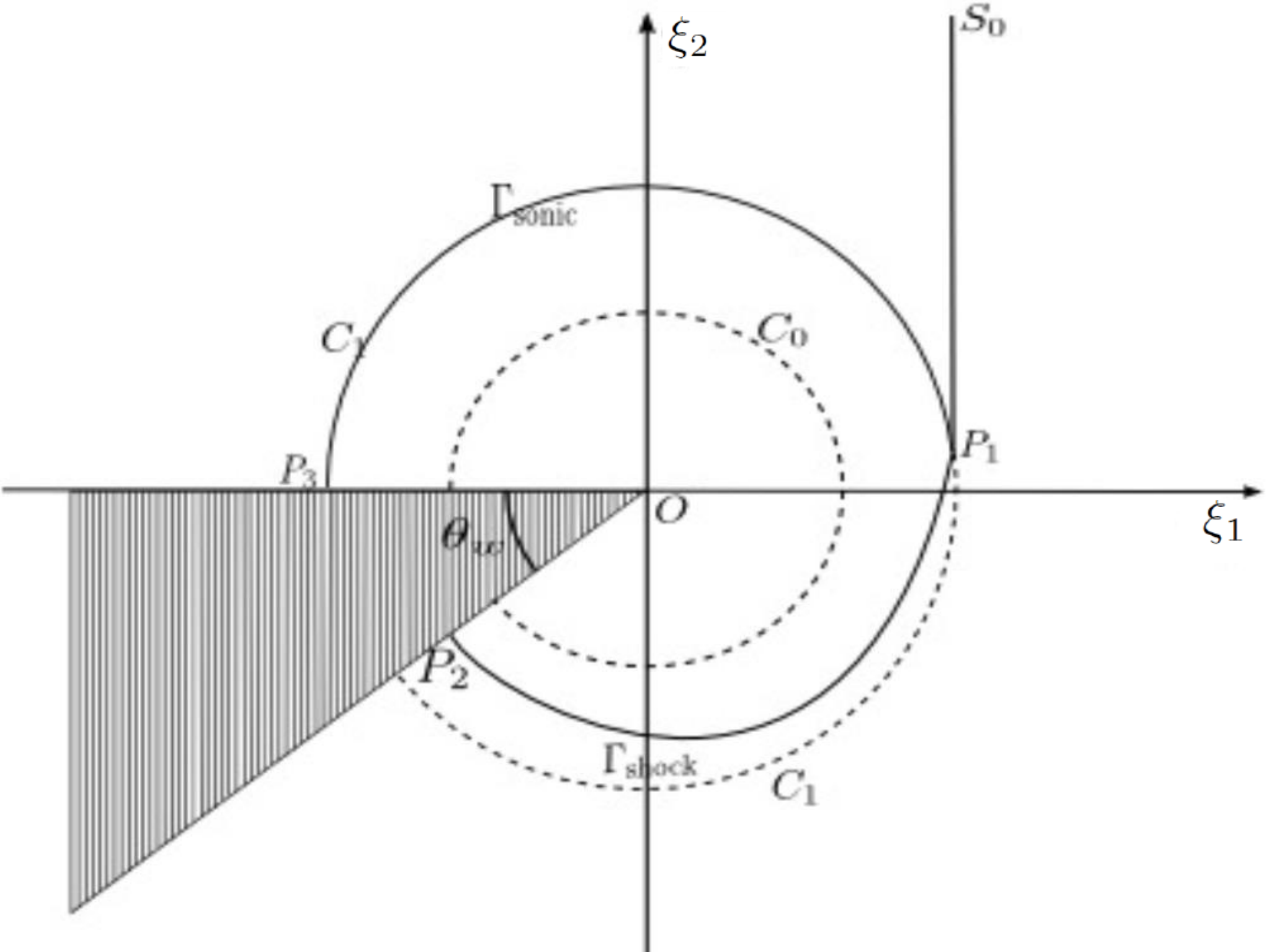}
  \caption{Shock diffraction configuration ({\it cf}. \cite{CDX})}
  \label{cog-2}
\end{figure}

Notice that Problem \ref{prob:4.1}
is invariant under the self-similar scaling:
$
(t, {\bf x})\to (\alpha t, \alpha {\bf x})$
for $\alpha\ne 0$.
In the self-similar $\xxi$--coordinates, system \eqref{1.1a}--\eqref{pressure}
can be rewritten in form \eqref{HCL-1b} with \eqref{1.1ab}.
In the polar coordinates $(r,\theta), r=\lvert\xxi\rvert$, the
system can be further written as
\begin{equation}\label{1.5b}
\partial_{r}\begin{pmatrix}
r\rho- m\cos\theta-n\sin\theta\\ rm-p(\rho)\cos\theta\\
rn-p(\rho)\sin\theta
 \end{pmatrix}+\partial_{\theta}
\begin{pmatrix}
m\sin\theta-n\cos\theta\\ p(\rho)\sin\theta\\
-p(\rho)\cos\theta
\end{pmatrix}
=\begin{pmatrix}
\rho+\frac{\cos\theta}{r}m+\frac{\sin\theta}{r}n \\  m+\frac{\cos\theta}{r}p(\rho)\\
n+\frac{\sin\theta}{r}p(\rho)
 \end{pmatrix}.
\end{equation}
The location of the incident shock $S_0$ for large $r\gg 1$ is:
\begin{equation}\label{ShockLocation:1}
\xi_1=\xi_1^0:=\sqrt{\frac{p(\rho_1)-p(\rho_0)}{\rho_1-\rho_0}}>0.
\end{equation}
Then Problem \ref{prob:4.1} can be reformulated as a boundary value problem in an
unbounded domain (see Fig. {\rm \ref{cog-2}}): $\,$
{\it Seek a solution of
system \eqref{HCL-1b} with \eqref{1.1ab},
or equivalently \eqref{1.5b}, with the
asymptotic boundary condition when $r\to \infty${\rm :}
\begin{equation}\label{initial-condition-2}
(\rho, m, n)\to
\begin{cases}
(\rho_0, 0, 0)\,
&\mbox{{\rm in} $\{\xi_1>\xi_1^0, \xi_2>0\}\cup\{-\pi+\theta_{\rm w}\le\arctan(\frac{\xi_2}{\xi_1})\le 0\}$},\\
(\rho_1, m_1, 0) \, &\mbox{{\rm in} $\{\xi_1<\xi_1^0,\xi_2>0\}$},
\end{cases}
\end{equation}
and the slip boundary condition along the wedge boundary $\partial
W${\rm :}
\begin{equation}\label{boundary-condition-ss}
(m, n)\cdot \nnu_{\rm w}\mid_{\partial W}=0.
\end{equation}
}

For a smooth solution $U=(\rho,m,n)$ of system \eqref{HCL-1b} with \eqref{1.1ab}, we may
eliminate $m$ and $n$ in \eqref{1.1a} to obtain a second-order
nonlinear equation for $\rho$: \begin{equation} \label{1.6a}
\big((c^{2}-\xi_1^{2})\rho_{\xi_1}-\xi_1\xi_2\rho_{\xi_2}+\xi_1\rho\big)_{\xi_1}
+\big((c^{2}-\xi_2^{2})\rho_{\xi_2}-\xi_1\xi_2\rho_{\xi_1}+\xi_2\rho\big)_{\xi_2}
-2\rho=0.
\end{equation}
Correspondingly, Eq. \eqref{1.6a} in the polar coordinates $(r,
\theta), r=\lvert\xxi\rvert$, takes the form
\begin{equation}\label{1.8b}
\big((c^{2}-r^{2})\rho_{r}\big)_{r}+\frac{c^{2}}{r}\rho_{r}
+\big(\frac{c^{2}}{r^{2}}\rho_{\theta}\big)_{\theta}=0.
\end{equation}

In the self-similar $\xxi$--coordinates, as the incident shock $S_0$ passes
through the wedge corner, $S_0$ interacts with the sonic circle
$\Gamma_{\rm sonic}$ of state (1): $r=r_1$, and becomes a transonic
diffracted shock $\Gamma_{\rm shock}$, and the flow in domain
$\Omega$ behind the shock and inside $\Gamma_{\rm sonic}$ becomes
subsonic.

Consider system \eqref{1.5b} in the polar coordinates. Then the
Rankine-Hugoniot relations, {\it i.e.}, the jump conditions, are
\begin{eqnarray*}
[p][\rho]=[m]^2+[n]^2, \qquad  \frac{{\rm d} r}{{\rm d}\theta}=
r\frac{\sqrt{r^2-\bar{c}^2(\rho,\rho_0)}}{\bar{c}(\rho,\rho_0)},
\end{eqnarray*}
with
$\bar{c}(\rho,\rho_0)=\sqrt{\frac{p(\rho)-p(\rho_0)}{\rho-\rho_0}}$,
where the plus branch has been chosen so that $\frac{dr}{d\theta}>0$.
Differentiating the first equation above along $\Gamma_{\text{shock}}$ and using
the equations obtained above, we have
\begin{equation}\label{2.41a}
\beta_1\rho_r+\beta_2\rho_{\theta}=0
\qquad \text{on $\Gamma_{\text{shock}}
:=\{(r(\theta),\theta)\, :\,\theta\in [\theta_{\rm w}, \theta_1]\}$}.
\end{equation}
where $\beta=(\beta_1,\beta_2)$ is a function of $(\rho_0, \rho,
r(\theta), r'(\theta))$ with
\begin{equation*}
\beta_{1}=r'(\theta)\big(c^{2}(r^2-\bar{c}^2)-3\bar{c}^2(c^{2}-r^2)\big),\quad
\beta_{2}=3c^{2}(r^{2}-\bar{c}^2)-\bar{c}^2(c^{2}-r^{2}).
\end{equation*}
Then the obliqueness becomes
\begin{equation*}
\mu:=\beta\cdot(1,-r'(\theta))=-2r^2(c^2-\bar{c}^{2})r'(\theta)
\ne 0,
\end{equation*}
where $(1,-r'(\theta))$ is the outward normal to $\Omega$ on
$\Gamma_{\text{shock}}$. Note that $\mu$ becomes zero when
$r'(\theta)=0$, {\it i.e.}, $r=\bar{c}(\rho,\rho_0)$, where
$$
\beta_1=0, \qquad \beta_2=-\bar{c}^2(c^2-r^2)<0,
$$
since $c^2(\rho)>\bar{c}^2(\rho,\rho_0)=r^2$ if $\rho>\rho_0$.

The second condition on $\Gamma_{\text{shock}}$ is the shock
equation: \begin{equation}\label{2.42a}
\frac{{\rm d} r}{{\rm d}\theta}=r\frac{\sqrt{r^2-\bar{c}^2(\rho,\rho_0)}}{\bar{c}(\rho,\rho_0)}
:=g(r,\theta,\rho(r,\theta)), \qquad\,\,\, r(\theta_1)=r_1,
\end{equation}
where $(r_1,\theta_1)$ are the polar coordinates of
$P_1=(\xi_1^0,\xi_2^0)$.

At point $P_2$, $r'(\theta_{\rm w})=0$, \eqref{2.41a} does not satisfy the oblique
derivative boundary condition. We may alternatively
express this as a one-point Dirichlet condition by solving
$r(\theta_{\rm w})=\bar{c}(\rho(r(\theta_{\rm w}), \theta_{\rm w}),\rho_0)$. In order
to deal with this equation, we use the  notation:
\begin{equation}
a=(\bar{c}_{b})^{-1}(r)\qquad  \text{when}\,\,
\bar{c}_b:=\bar{c}(a,b)=r \,\,\,\mbox{for fixed $b$},
\end{equation}
so that
\begin{equation} \label{2.47a}
\rho(P_1)=\bar{\rho}=(\bar{c}_{\rho_0})^{-1}(r(\theta_{\rm w})).
\end{equation}

The boundary condition on the wedge is the slip boundary condition,
{\it i.e.}, $(m,n)\cdot\nnu_{\rm w}=0$. Differentiating it along the wedge
and combining this with the second and third equations in
\eqref{1.1a}, we conclude that $\rho$ satisfies
\begin{equation}
\label{2.45a} \rho_{\nnu_{\rm w}}=0\qquad
\text{on $\Gamma_0:=\partial\Omega\cap(\{\theta=\pi\}\cup \{\theta=\theta_{\rm w}\})$}.
\end{equation}

The Dirichlet boundary condition on $\Gamma_{\rm sonic}$ is:
\begin{equation} \label{2.43a}
\rho=\rho_1\qquad
\text{on $\Gamma_{\text{sonic}}:=\partial\Omega\cap\partial B_{c_1}(0)$}.
\end{equation}
On the Dirichlet boundary $\Gamma_{\rm sonic}$, Eq. \eqref{1.8b}
becomes degenerate elliptic from the inside of
$\Omega$.

With the derivation of the free boundary conditions on
$\Gamma_{\rm shock}$ and the fixed boundary conditions on
$\Gamma_{\rm sonic}\cup \Gamma_0$,
Problem \ref{prob:4.1} is further
reduced to the following free boundary problem for Eq. \eqref{1.8b}
in domain $\Omega$, with $(m, n)$
correspondingly determined by \eqref{1.5b}.

\medskip
\begin{problem}[Free Boundary Problem]\label{prob:4.2}
{\it Seek a solution $(\rho(r, \theta), r(\theta))$ such that
$\rho(r,\theta)$ and $r(\theta)$ are determined by
Eq. \eqref{1.8b} in domain $\Omega$ and the free boundary
conditions \eqref{2.41a}--\eqref{2.47a}
on $\Gamma_{\rm shock}=\{(r(\theta),\theta)\,:\, \theta_{\rm w}\le \theta\le \theta_1\}$,
in addition to the
Neumann boundary condition \eqref{2.45a} on wedge $\Gamma_0$
and the Dirichlet boundary condition \eqref{2.43a} on the degenerate
boundary $\Gamma_{\rm sonic}$, the sonic circle of state {\rm (1)} $(${\it
cf.} Fig. {\rm \ref{cog-2}}$)$.}
\end{problem}

\subsection{Global Solutions of Riemann Problem II:
\\Free Boundary Problem, Problem \ref{prob:4.2}}

To solve Riemann Problem II, it suffices to deal with the free boundary problem,
Problem \ref{prob:4.2}, which has been solved as stated in the following theorem.

\begin{theorem}[Chen-Deng-Xiang \cite{CDX}]\label{1}\label{2.1}
Let the wedge angle $\theta_{\rm w}$ be between $-\pi$ and $0$. Then there
exists a global solution, a density function $\rho(r,\theta)$ in domain $\Omega$,
and a free boundary $\Gamma_{\rm shock}=\{(r(\theta),\theta)\,:\, \theta_{\rm w}\le\theta\le \theta_1\}$,
of {\rm Problem \ref{prob:4.2}} such that
$$
\rho\in C^{2+\alpha}(\Omega)\cap C^{\alpha}(\overline{\Omega}),
\quad r\in C^{2+\alpha}([\theta_{\rm w},\theta_1))\cap
C^{1,1}([\theta_{\rm w},\theta_1]).
$$
Moreover, solution $(\rho(r,\theta), r(\theta))$ satisfies the
following properties{\rm :}

\begin{enumerate}
\item[\rm (i)] $\rho> \rho_0$ on the free boundary $\Gamma_{\rm shock}${\rm ;} that is, $\Gamma_{\rm shock}$
is separated from the sonic circle $C_0$ of state {\rm (0)}.

\item[\rm (ii)] The free boundary $\Gamma_{\rm shock}$ is strictly convex up to point $P_1$,
except point $P_2$, in the self-similar $\xxi$--coordinates.

\item[\rm (iii)] The density function $\rho(r,\theta)$ is $C^{1,\alpha}$ up to $\Gamma_{\rm sonic}$ and
Lipschitz continuous across $\Gamma_{\rm sonic}$.

\item[\rm (iv)] The Lipschitz regularity of $\rho(r,\theta)$
across $\Gamma_{\rm sonic}$ and at $P_1$ from the inside is optimal.
\end{enumerate}
\end{theorem}

\smallskip
Similar to the proof of Theorem \ref{thm1},
Theorem \ref{1} is established in two steps. First,
the regularized approximate free boundary problem for \eqref{1.8b} involving
two small parameters $\varepsilon$ and $\delta$ is solved.
Then the limits: $\varepsilon\rightarrow0$ and
$\delta\rightarrow0$ are proved to yield a solution of
Problem \ref{prob:4.2}, {\it i.e.}, \eqref{1.8b}--\eqref{2.43a},
with the optimal regularity.

In Theorem \ref{1}, a global
solution $\rho$ of Eq. \eqref{1.8b}
in $\Omega$
is constructed, by combining this function with $\rho=\rho_1$ in state (1)
and $\rho=\rho_0$ in state (0). That is, the global
density function $\rho$ that is piecewise constant in the supersonic
domain is Lipschitz continuous across the degenerate sonic
boundary $\Gamma_{\rm sonic}$ from $\Omega$ to state (1).
To recover the momentum vector function $(m,n)$, we can integrate
the second and third equations in \eqref{1.5b}. These can also be
written in the radial variable $r$,
\begin{equation} \label{6.1}
\frac{\partial (m,n)}{\partial r}=\frac{1}{r} \D p(\rho)
\end{equation}
and integrated from the boundary of the subsonic domain
toward the origin.

It has been proved that the limit of $\D\rho$ does not exist at
$P_1$ as $\xxi$ in $\Omega$ tends to $\xxi^0$, but $\lvert \D
c(\rho)\rvert$ has a upper bound. Thus, $p(\rho)$ is Lipschitz, which
implies that $(m, n)$ is at least Lipschitz across the sonic circle
$\Gamma_{\rm sonic}$.
Furthermore, $(m,n)$ has the same regularity as $\rho$ inside
$\Omega$, except for origin $r=0$. However, $(m, n)$ may be
multi-valued at origin $r=0$. Therefore, we have

\begin{theorem}[Chen-Deng-Xiang \cite{CDX}]\label{6.1a}
Let the wedge angle $\theta_{\rm w}$ be between $-\pi$ and $0$. Then there
exists a global solution $(\rho, m, n)(r,\theta)$ with the diffracted shock
 $\Gamma_{\rm shock}=\{(r(\theta),\theta)\,:\, \theta_{\rm w}\le \theta\le \theta_1\}$ of {\rm
Problem \ref{prob:4.2}} such that
$$
(\rho, m, n)\in C^{2+\alpha}(\Omega), \quad \rho\in
C^{\alpha}(\overline{\Omega}), \quad r\in
C^{2+\alpha}([\theta_{\rm w},\theta_1))\cap C^{1,1}([\theta_{\rm w},\theta_1]),
$$
and $(\rho, m, n)=(\rho_1, m_1, 0)$ in domain $\{\xi_1<\xi_1^0,\,
r>r_1\}$ and $(\rho_0, 0,0)$ in domain $\{\xi_1>\xi_1^0,\,
\xi_2>\xi_2^0\}\cup\{r>r(\theta),\,\theta_{\rm w}\le\theta\le \theta_1\}$.
Moreover, solution $(\rho,m,n)(r,\theta)$ with the diffracted shock
$\Gamma_{\rm shock}$
satisfies properties {\rm (i)}--{\rm (ii)} in {\rm Theorem \ref{2.1}} and

\smallskip
\begin{enumerate}
\item[\rm (i)] $(\rho,m,n)$ is $C^{1,\alpha}$ up to $\Gamma_{\rm sonic}$ and
Lipschitz continuous across $\Gamma_{\text{\rm sonic}}$.

\item[\rm (ii)] The Lipschitz regularity of solution $(\rho,m,n)$
across $\Gamma_{\rm sonic}$ and at $P_1$ from the inside is optimal.

\item[\rm (iii)] The momentum vector function $(m,n)$ may be multi-valued at the origin.
\end{enumerate}
\end{theorem}

In particular, Theorem \ref{6.1a} implies the following facts:

\smallskip
\begin{enumerate}
\item[(a)] The optimal regularity of $(\rho, m, n)(r, \theta)$
across $\Gamma_{\rm sonic}$ and at $P_1$ from the inside is $C^{0,1}$,
{\it i.e.}, Lipschitz continuity.

\item[(b)] The diffracted shock $\Gamma_{\rm shock}$ is definitely not degenerate at point $P_2$.
This had been an open question even when the wedge angle is
$\frac{\pi}{2}$ as in \cite{e}, though it is physically plausible.

\item[(c)] The
diffracted shock $\Gamma_{\rm shock}$
away from point $P_2$ is strictly convex
and has a jump at point $P_1$ from a positive value to zero, while
the strict convexity of $\Gamma_{\rm shock}$
fails at $P_2$.
\end{enumerate}

More details can be found in Chen-Deng-Xiang \cite{CDX}. Similar results
can be obtained for the pressure gradient equation introduced in Section 3 above.
In Chen-Feldman-Hu-Wang \cite{ChenFeldmanHuXiang}, the loss of regularity of solutions
of Problem \ref{prob:4.1} for  the potential flow equation \eqref{1-a}--\eqref{1-c}, or \eqref{HCL-2}
with \eqref{1-b1}, has been shown, which implies that
the solution configuration for this case is much more complicated.

\section{Two-Dimensional Riemann Problem III: \\
The Prandtl-Meyer Problem for Unsteady
Supersonic Flow onto Solid Wedges \\ for the Euler Equations for Potential Flow}
\label{PrandtlMeyerProblSect}

Now we present the third Riemann problem, Riemann Problem III, for the Prandtl-Meyer problem for unsteady
supersonic flow onto solid wedges for the Euler equations for potential flow
in form \eqref{HCL-2} with \eqref{1-b1}--\eqref{HCL-2e}, or \eqref{1-a}--\eqref{1-c};
see also \cite{BCF-14,EllingLiu1,Meyer,Prandtl}.

\subsection{2-D Riemann Problem III: The Prandtl-Meyer Problem for Unsteady Supersonic Flow onto \\ Solid Wedges for Potential Flow}

Consider a supersonic flow
with the constant density $\irho >0$ and velocity ${\bf u}_{0}=(u_0, 0)$,
$u_0> c_0:=c(\rho_0)$, which
impinges toward a symmetric wedge:
\begin{equation}\label{wedge-symm}
W:= \big\{(x_1,x_2)\,:\, \lvert x_2\rvert < x_1 \tan \theta_{\rm w}, x_1 > 0\big\}
\end{equation}
at $t = 0$.
If $\theta_{\rm w}$ is less
than the detachment angle $\theta_{\rm w}^{\rm d}$,
then the well-known {\emph{shock polar analysis}} demonstrates that there are two different
steady weak solutions:
{\emph{the steady solution}} $\bar{\Phi}$ of {\emph{weaker shock strength}} and {\emph{the steady solution}} of {\emph{stronger shock strength}},
both of which satisfy the entropy condition
and the slip boundary condition (see Fig.  \ref{Figure2}); see also \cite{BCF-14,Chen2,CFr}.
Then the dynamic stability of the steady transonic solution  $\bar{\Phi}$  of weaker shock strength for potential flow can be
formulated as
the following problem:

\begin{figure}
 	\centering
 	\includegraphics[height=32mm]{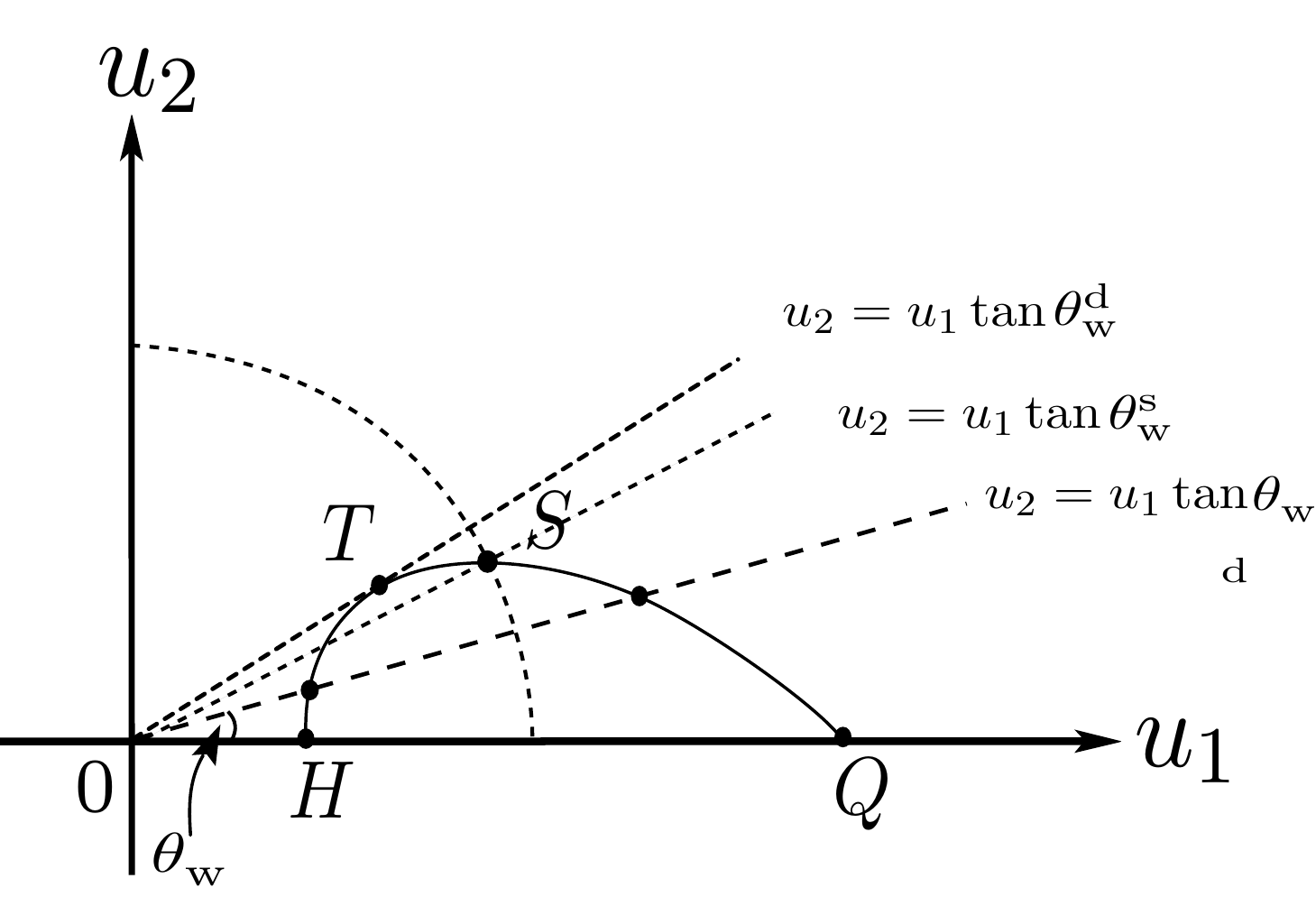}
\hspace{6mm}
 		\includegraphics[height=31mm]{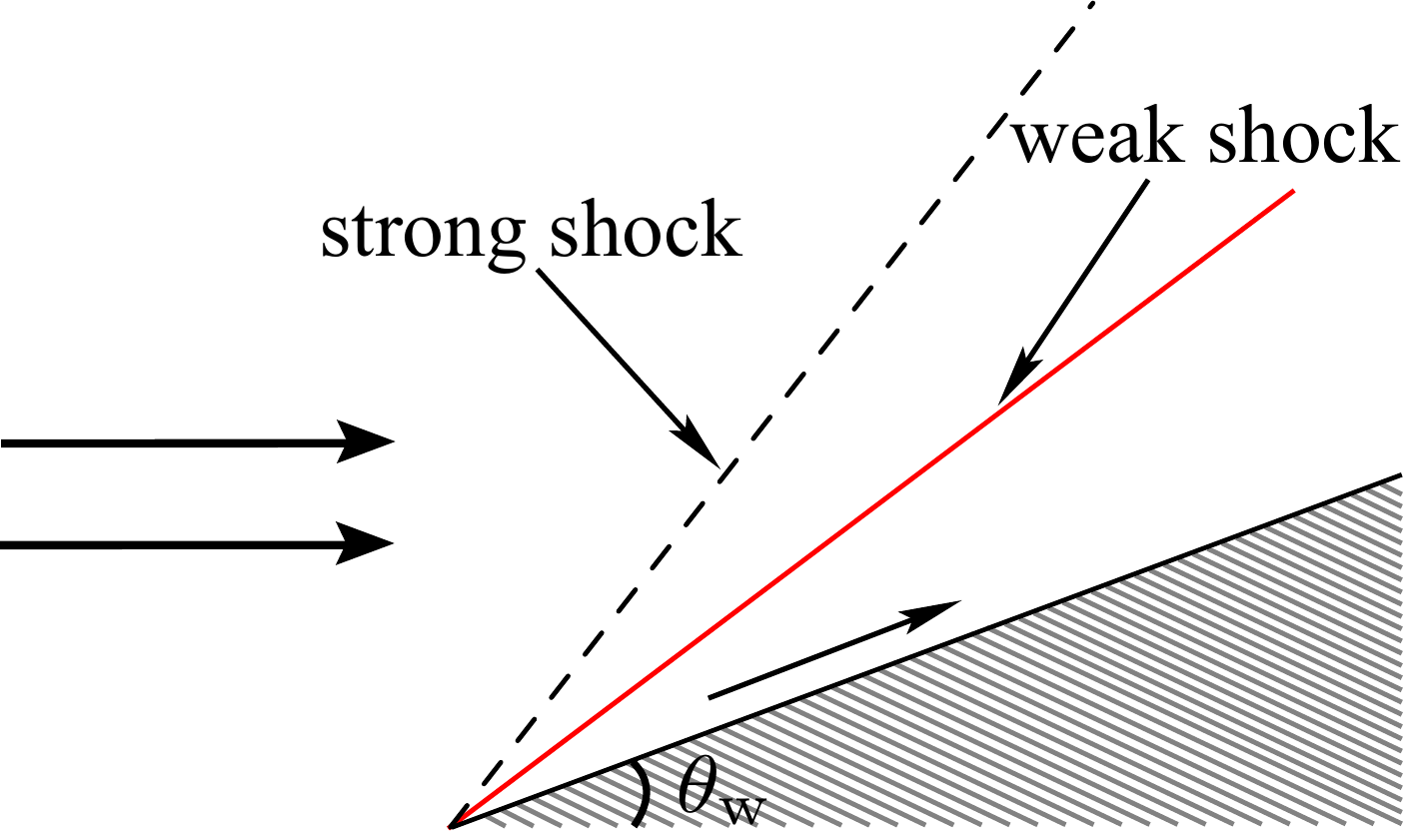}
		\caption{The shock polar in the $\uu$-plane and uniform steady (weak/strong) shock flows (see \cite{Chen2})}
		\label{Figure2}
 \end{figure}

\begin{problem}[Riemann Problem III:  The Prandtl-Meyer Problem for Unsteady
Supersonic Flow onto Solid Wedges]\label{problem-1}
Given $\gam>1$, fix $(\irho, u_0)$ with $u_0>c_0$.
For a fixed $\theta_{\rm w}\in (0,\theta_{\rm w}^{\rm d})$,
seek a global entropy solution $\Phi\in W^{1,\infty}_{\rm loc}(\RR_+\times (\RR^2\setminus W))$
of Eq. \eqref{HCL-2} with \eqref{1-b1}--\eqref{HCL-2e}
and $B=\frac{u_0^2}{2}+h(\irho)$ so that $\Phi$ satisfies the initial condition at $t=0${\rm :}
\begin{equation}\label{1-d}
(\rho,\Phi)\rvert_{t=0}=(\irho, u_0 x_1) \qquad\,\,\, \text{for}\;\;\xx\in \RR^2\setminus W,
\end{equation}
and the slip boundary condition along the wedge boundary $\der W${\rm :}
\begin{equation}\label{1-e}
\nabla_{\bf x}\Phi\cdot \nnu_{\rm w} \rvert_{\der W}=0,
\end{equation}
where $\nnu_{\rm w}$ is the exterior unit normal to $\der W$.
In particular, we seek a solution $\Phi\in W^{1,\infty}_{\rm loc}(\RR_+\times (\RR^2\setminus W))$
that converges to
the steady  solution $\bPhi$ of weaker oblique shock strength
corresponding to
the fixed parameters $(\irho, u_0, \gam, \theta_{\rm w})$
with $\bar{\rho}=h^{-1}(B-\frac 12\lvert \nabla \bPhi\rvert^2)$,
when $t\to \infty$, in the following sense{\rm :}
For any $R>0$, $\Phi$ satisfies
\begin{equation}
\label{time-asymp-lmt}
\lim_{t\to \infty} \|(\nabla_{\bf x}\Phi(t,\cdot)-\nabla_{\bf x}\bPhi,
\rho(t,\cdot)-\bar{\rho})\|_{L^1(B_R({\bf 0})\setminus W)}=0
\end{equation}
for $\rho(t,{\bf x})$ given by \eqref{1-b1}.
\end{problem}

\smallskip
Since the initial data in \eqref{1-d} do not satisfy the boundary condition \eqref{1-e},
a boundary layer is generated along the wedge boundary starting at $t=0$,
which forms the Prandtl-Meyer reflection configurations; see
Bae-Chen-Feldman \cite{BCF-14} and the references cited therein.

In order to define the notion of weak solutions of Problem \ref{problem-1}, it is noted that
the boundary condition can be written as $\rho\nabla_{\bf x}\Phi\cdot\nnu_{\rm w}=0$ on $\partial W$,
which is the spatial conormal condition for Eq. \eqref{HCL-2} with \eqref{1-b1}--\eqref{HCL-2e}.
Then we have

\begin{definition}[Weak Solutions of {Problem \ref{problem-1}}: Riemann Problem III]\label{weakSol-def-Prob1}
A function
$
\Phi\in W^{1,1}_{\rm loc}(\RR_+\times (\bR^2\setminus W))
$
is called a weak
solution of  {\rm  Problem \ref{problem-1}}
if $\Phi$ satisfies the following properties{\rm :}

\smallskip
\begin{enumerate}[\rm (i)]
\item
\label{weakSol-def-i1a-Prob1}
$B-\big(\partial_t\Phi+\frac{1}{2}\lvert \nabla_\xx\Phi\rvert ^2\big)\ge h(0+)\,\,$ {\it a.e.} in
$\RR_+\times (\RR^2\setminus W)$,

\smallskip
\item
\label{weakSol-def-i2a-Prob1}
For $\rho(\partial_t \Phi, \nabla_\xx\Phi)$ determined by \eqref{1-b1},
$$
(\rho(\partial_t\Phi, \lvert \nabla_\xx\Phi\rvert ^2), \rho(\partial_t\Phi, \lvert \nabla_\xx\Phi\rvert ^2)\lvert \nabla_\xx\Phi\rvert )\in
(L^1_{\rm loc}(\RR_+\times \overline{\RR^2\setminus W}))^2,
$$

\item
\label{weakSol-def-i3a-Prob1}
For every $\zeta\in C^\infty_c({\RR_+}\times \bR^2)$,
\begin{align*}
&\int_0^\infty\int_{\RR^2\setminus W}\Big(\rho(\partial_t\Phi, \lvert \nabla_\xx\Phi\rvert ^2)\partial_t\zeta
 +\rho(\partial_t\Phi, \lvert \nabla_\xx\Phi\rvert ^2)
\nabla_\xx\Phi\cdot\nabla_\xx\zeta\Big) \,\dr \xx \dr t\\
& +\int_{\RR^2\setminus W}\rho_0\zeta(0, \xx)\,\dr\xx=0.
\end{align*}
\end{enumerate}
\end{definition}

Since $\zeta$ does not need to be zero on $\partial\Lambda$,
the integral identity in Definition {\rm \ref{weakSol-def-Prob1}}
is a weak form of equation \eqref{HCL-2} with \eqref{1-b1}--\eqref{HCL-2e}
and the boundary condition $\rho\nabla_{\bf x}\Phi\cdot\nnu_{\rm w}=0$ on $\partial W$.
A weak solution is called an entropy solution if it satisfies the entropy condition
that is consistent with the second law of thermodynamics
$(${\it cf}. {\rm \cite{CF-book2018,CFr,Da,Lax}}$)$. In particular, a piecewise smooth solution
is an entropy solution if the discontinuities are all shocks.

\subsection{Reformulation of Riemann Problem III}

Notice that Eq. \eqref{HCL-2} with  \eqref{1-b1}--\eqref{HCL-2e}
is invariant under the self-similar scaling \eqref{sss},
so that it admits  self-similar solutions in form \eqref{ssform}.
Then the pseudo-potential function $\varphi=\phi-\frac{1}{2}\lvert \xxi\rvert^2$
satisfies
the following equation:
\begin{equation}
\label{2-1}
{\rm div}(\rho(\lvert \D\vphi\rvert ^2,\vphi)\D\vphi)+2\rho(\lvert \D\vphi\rvert^2,\vphi)=0
\end{equation}
for
\begin{equation}
\label{1-o}
\rho(\lvert \D\vphi\rvert^2,\vphi)=
\bigl(B_0- (\gamma-1)(\frac{1}{2}\lvert \D\vphi\rvert^2+\varphi)\bigr)^{\frac{1}{\gam-1}},
\end{equation}
where $B_0=(\gamma-1)B+1$.
Eq. \eqref{2-1}
written in the non-divergence form is
\begin{equation}\label{nondivMainEq}
(c^2-\varphi_{\cxi}^2)\varphi_{\cxi\cxi}-2\varphi_\cxi\varphi_\ceta\varphi_{\cxi\ceta}
+(c^2-\varphi_{\ceta}^2)\varphi_{\ceta\ceta}+2c^2-\lvert \D\varphi\rvert^2=0,
\end{equation}
where the sonic speed $c=c(\lvert \D\varphi\rvert ^2,\varphi)$ is determined by
\begin{equation}\label{c-through-density-function}
c^2(\lvert \D\varphi\rvert^2,\varphi)=
\rho^{\gamma-1}(\lvert \D\varphi\rvert^2,\varphi)
=B_0-(\gamma-1)\big(\frac{1}{2}\lvert \D\varphi\rvert^2+\varphi\big).
\end{equation}

Eq. \eqref{2-1} is a nonlinear PDE of mixed elliptic-hyperbolic type. It is elliptic
at $\xxi$ if and only if
\begin{equation}
\label{1-f}
\lvert \D\vphi\rvert<c(\lvert \D\vphi\rvert^2,\vphi)   \qquad \,\, \mbox{at $\xxi$},
\end{equation}
and is hyperbolic if the opposite inequality holds.

One class of solutions of (\ref{2-1})
is that of {\em constant states} which are the solutions
with constant velocity $\mathbf{v}\in \mathbb{R}^2$.
Then the pseudo-potential of the constant state $\mathbf{v}$ satisfies
$\D\varphi=\mathbf{v}-\xxi$ so that
\begin{equation}\label{constantStatesForm}
\varphi(\xxi)=-\frac 12\lvert \xxi\rvert^2+\mathbf{v}\cdot\xxi +C,
\end{equation}
where $C$ is a constant. For such $\varphi$, the expressions
in \eqref{1-o} and \eqref{c-through-density-function}
imply that the density and sonic
speed are positive constants $\rho$ and $c$, {\it i.e.},
independent of $\xxi$.
Then, from \eqref{1-f}--\eqref{constantStatesForm},
the ellipticity condition for the constant state $\mathbf{v}$ is
$$
\lvert \xxi -\mathbf{v}\rvert<c.
$$
Thus,
Eq. (\ref{2-1})
is elliptic inside the {\em sonic circle}
with center $\mathbf{v}$ and radius $c$, and hyperbolic outside this circle.

Moreover,
if density $\rho$ is a constant, then the solution is also a constant state;
that is, the corresponding
pseudo-potential $\vphi$ is of form \eqref{constantStatesForm}.

\smallskip
Since the problem involves transonic shocks, we have to consider
weak solutions of Eq. \eqref{2-1},
which admit shocks.
A shock is a curve across which $\D\vphi$ is discontinuous. If $\Lambda^+$ and $\Lambda^-(:=\Lambda\setminus \ol{\Lambda^+})$
are two nonempty open subsets of a domain $\Lambda\subset \RR^2$, and ${S}:=\der\Lambda^+\cap \Lambda$ is a $C^1$-curve
across which $\D\vphi$ has a jump, then $\vphi\in W^{1,1}_{\rm loc}\cap C^1(\Lambda^{\pm}\cup S)\cap C^2(\Lambda^{\pm})$
is a global weak solution of \eqref{2-1} in $\Lambda$ if and only if $\vphi$ is in $W^{1,\infty}_{\rm loc}(\Lambda)$
and satisfies Eq. \eqref{2-1} and the Rankine-Hugoniot conditions on $S$:
\begin{align}
&\vphi\rvert_{\Lambda^+\cap S}=\vphi\rvert_{\Lambda^-\cap S}, \label{1-i}\\
&\rho(\lvert \D\vphi\rvert^2, \vphi)\D\vphi\cdot\nnu_{\rm s}\rvert_{\Lambda^+\cap S}
=\rho(\lvert \D\vphi\rvert^2, \vphi)\D\vphi\cdot\nnu_{\rm s}\rvert_{\Lambda^-\cap S}. \label{1-h}
\end{align}
A piecewise smooth solution with the discontinuities is called an {\it entropy solution}
of \eqref{2-1} if it satisfies the entropy condition: {\it density $\rho$ increases in the pseudo-flow direction of
$\D\varphi\rvert_{\Lambda^+\cap S}$ across the discontinuity}. Then such a discontinuity is called a shock.

\medskip
As the upstream flow has the constant velocity $\uu_0=(u_0,0)$,
the corresponding pseudo-potential
$\ivphi$ has the expression of
\begin{equation}
\label{1-m}
\ivphi=-\frac 12\lvert \xxi\rvert^2+u_0\xi_1
\end{equation}
directly from \eqref{constantStatesForm} with the choice of $B$ in Problem \ref{problem-1}.
Since
the symmetry of the domain and the upstream flow in {Problem \ref{problem-1}}
with respect to the $x_1$--axis,
{Problem \ref{problem-1}} can then be reformulated as the following boundary value problem
in the domain:
$$
\Lambda:=\RR^2_+\setminus\{\xxi\,: \,\xi_2\le \xi_1\tan\theta_{\rm w},\, \xi_1\ge 0\}
$$
in the
self-similar coordinates $\xxi$,
which corresponds to domain $\{(t, {\bf x})\, :\, {\bf x}\in \RR^2_+\setminus W,\, t>0\}$
in the $(t, {\bf x})$--coordinates, where $\RR^2_+=\{\xxi\,: \,\xi_2>0\}$:
{\it $\,$ Seek a solution $\vphi$ of Eq. \eqref{2-1} in the self-similar domain $\Lambda$
with the slip boundary condition{\rm :}
\begin{equation}
\label{1-k}
\D\vphi\cdot \nnu_{\rm w}\rvert_{\partial \Lambda}=0
\end{equation}
and the asymptotic boundary condition{\rm :}
\begin{equation}\label{1-k-b}
\vphi-\vphi_0\longrightarrow 0
\end{equation}
along each ray $R_\theta:=\{ \xi_1=\xi_2 \cot \theta, \xi_2 > 0 \}$
with $\theta\in (\theta_{\rm w}, \pi)$ as $\xi_2\to \infty$
in the sense that
\begin{equation}\label{1-k-c}
\lim_{r\to \infty} \|\varphi  - \varphi_0\|_{C(R_{\theta}\setminus B_r(0))} = 0.
\end{equation}
}

\begin{figure}
 	\centering
 	\includegraphics[height=32mm]{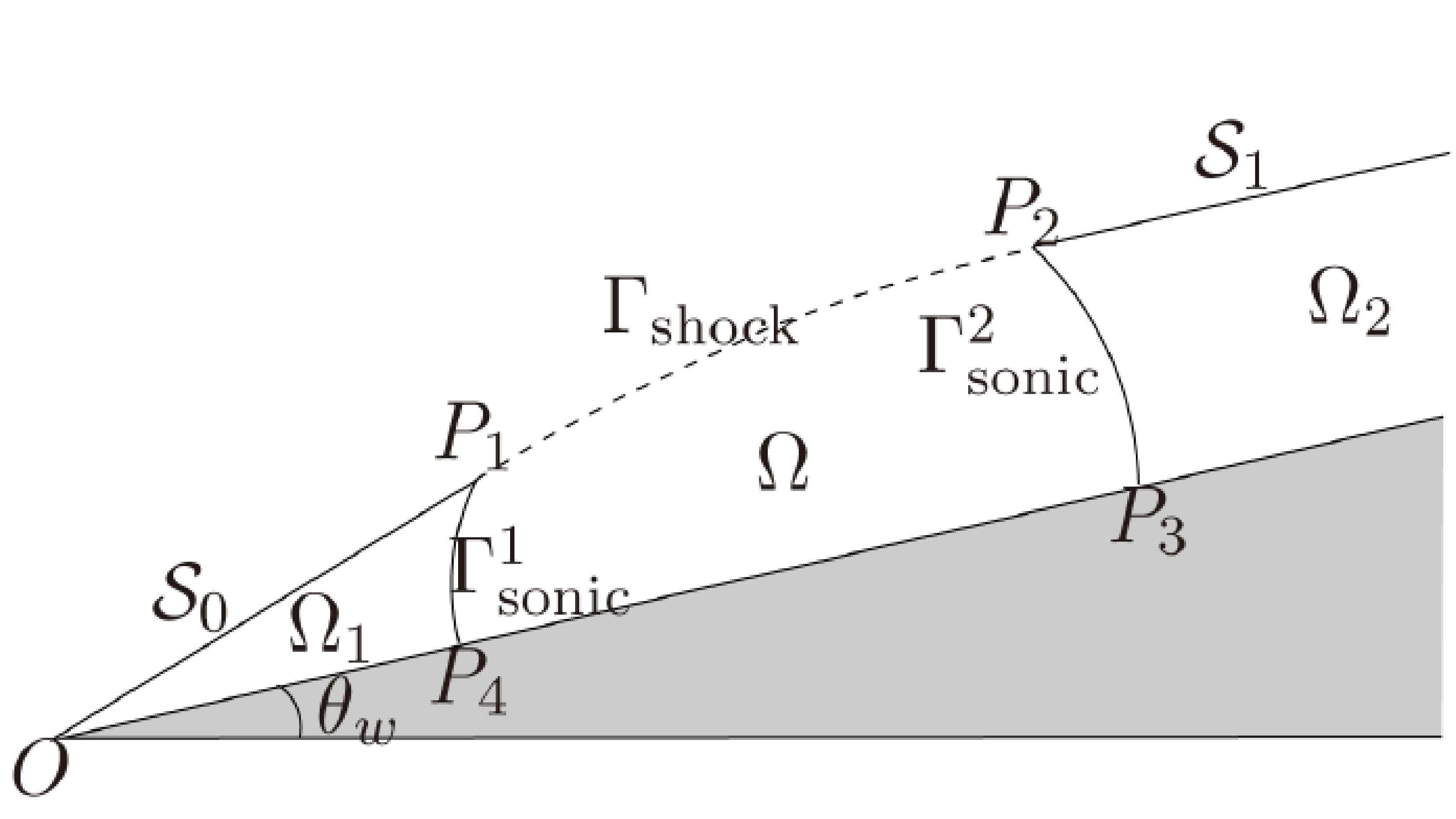}
	\caption{Self-similar solutions for $\theta_{\rm w}\in (0, \theta_{\rm w}^{\rm s})$
in the self-similar coordinates $\xxi$ ({\it cf}.  \cite{BCF-14})}\label{fig:global-structure-1}
 \end{figure}

\begin{figure}
\centering
\includegraphics[height=32mm]{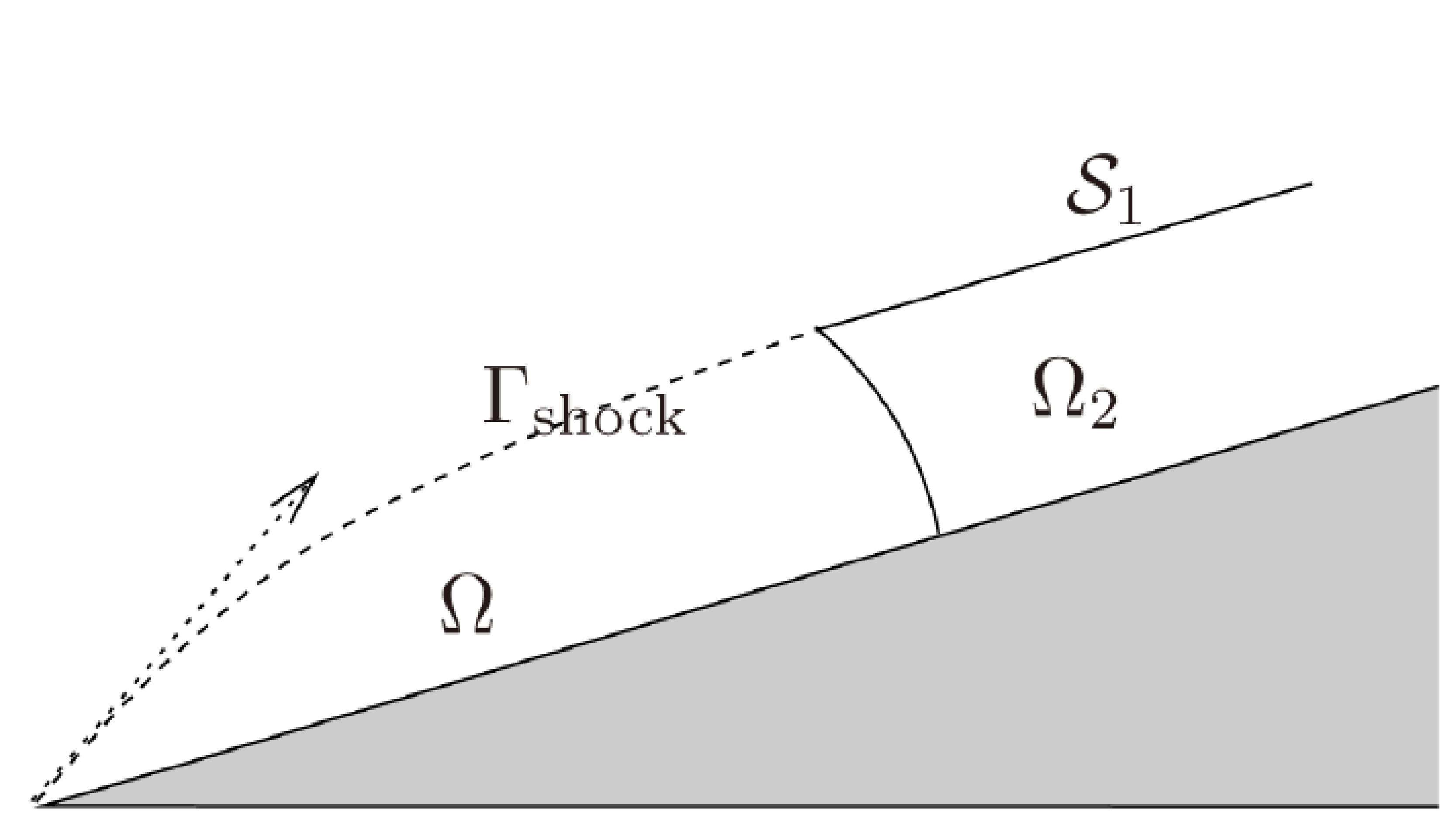}
\caption{Self-similar solutions for $\theta_{\rm w}\in [\theta_{\rm w}^{\rm s},\theta_{\rm w}^{\rm d})$
in the self-similar coordinates $\xxi$ ({\it cf}. \cite{BCF-14})}\label{fig:global-structure-2}
\end{figure}

Given $M_0>1$, $\rho_1$ and $\uu_1$ are determined via the shock polar as shown in Fig. \ref{Figure2}
for steady potential flow.
For any wedge angle $\theta_{\rm w}\in (0,\theta_{\rm w}^{\rm s})$,
line $v=u\tan\theta_{\rm w}$ and the shock polar
intersect at a point $\uu_1=(u_1,v_1)$ with $\lvert \uu_1\rvert >c_1$ and $u_{1}<u_0$;
while, for any $\theta_{\rm w}\in [\theta_{\rm w}^{\rm s}, \theta_{\rm w}^{\rm d})$,
they intersect at a point $\uu_1$ with $u_{1}>u_{\rm d}$
and $\lvert \uu_1\rvert<c_1$
where $u_{\rm d}$ is the $u$--component
of the unique detachment state $\uu_{\rm d}$ when $\theta_{\rm w}=\theta_{\rm w}^{\rm d}$.
The intersection state $\uu_1$ is the velocity for steady potential flow
behind an oblique shock ${S}_0$ attached to the wedge vertex with angle $\theta_{\rm w}$.
The strength of shock ${S}_0$ is relatively weak compared to the shock
given by the other intersection point on the shock polar,
hence ${S}_0$ is called
a \emph{weak oblique shock}
and the corresponding state $\uu_1$ is a \emph{weak state}.
Moreover, such a state
$\uu_1$ depends smoothly on $(u_0, \theta_{\rm w})$ and
is supersonic when $\theta_{\rm w}\in (0,\theta_{\rm w}^{\rm s})$
and subsonic when $\theta_{\rm w}\in [\theta_{\rm w}^{\rm s}, \theta_{\rm w}^{\rm d})$.

Once $\uu_1$ is determined, by \eqref{1-i}--\eqref{1-m},
the pseudo-potential $\vphi_1$ below the weak oblique shock ${S}_0$
is
\begin{equation}\label{1-n}
\vphi_1=-\frac 12\lvert \xxi\rvert^2+ \uu_1\cdot\xxi.
\end{equation}

We seek a global entropy solution
with two types of Prandtl-Meyer reflection configurations whose occurrence
is determined by the wedge angle $\theta_{\rm w}$ for the two different cases:
One contains a straight weak oblique shock ${S}_0$ attached to the wedge vertex $O$ and
connected to a normal shock ${S}_1$ through a curved shock  $\Gamma_{\rm shock}$
when $\theta_{\rm w}<\theta_{\rm w}^{\rm s}$,
as shown in Fig. \ref{fig:global-structure-1};
the other contains a curved shock  $\Gamma_{\rm shock}$ attached
to the wedge vertex
and connected to a normal shock ${S}_1$
when $\theta_{\rm w}^{\rm s}\le \theta_{\rm w}<\theta_{\rm w}^{\rm d}$,
as shown in Fig. \ref{fig:global-structure-2},
in which the curved shock $\Gamma_{\rm shock}$ is tangential to
the straight weak oblique shock  ${S_0}$ at the wedge vertex. To achieve these, we need
to compute the pseudo-potential function $\vphi$ below ${S_0}$.

By \eqref{1-i}--\eqref{1-k},
the pseudo-potential $\vphi_2$ below the normal shock ${S}_1$
is of the form:
\begin{equation}\label{1-n-b}
\vphi_2=-\frac 12\lvert \xxi\rvert^2+ \uu_2\cdot \xxi+ k_2
\end{equation}
for constant state $\uu_2$ and constant $k_2$; see \eqref{constantStatesForm}.
Then it follows
from \eqref{1-o} and \eqref{1-n}--\eqref{1-n-b} that the corresponding densities $\rho_1$ and $\rho_2$
are constants in the form:
\begin{equation}\label{2-n1}
\rho_k^{\gam-1}=\rho_0^{\gamma-1}+\frac{\gam-1}{2}\big(u_0^2-\lvert \uu_k\rvert^2\big)\qquad\,\,\mbox{for $k=1,2$}.
\end{equation}

Denote $\Wedge :=\partial W\cap\partial\Lambda$,
and the sonic arcs $\Gamma_{\rm sonic}^1:=P_1P_4$
on Fig. {\rm \ref{fig:global-structure-1}} and $\Gamma_{\rm sonic}^2:=P_2P_3$ on
Figs. {\rm \ref{fig:global-structure-1}}--{\rm \ref{fig:global-structure-2}}.
The sonic circle $\partial B_{c_1}(\uu_1)$ of the uniform state $\varphi_1$ intersects line ${S_0}$,
where $c_1=\rho_1^{\frac{\gamma-1}2}$ by \eqref{c-through-density-function}.
For the supersonic case $\theta_{\rm w}\in (0,\theta_{\rm w}^{\rm s})$,
there are two arcs of this sonic circle between ${S_0}$ and $\Wedge$ in $\Lambda$.
Note that $\Gamma_{\rm sonic}^1$ tends to point $O$ as $\theta_{\rm w} \nearrow \theta_{\rm w}^{\rm s}$
and is outside of $\Lambda$ for the subsonic case $\theta_{\rm w}\in [\theta_{\rm w}^{\rm s}, \theta_{\rm w}^{\rm d})$.
Similarly, the sonic circle $\partial B_{c_2}(\uu_2)$ of the uniform state $\varphi_2$ intersects line ${S_1}$,
where $c_2=\rho_2^{\frac{\gamma-1}2}$. There are two
arcs of this circle between ${S_1}$ and the line containing $\Wedge$.
Notice that $\varphi_1>\varphi_2$ on $\overline{\Gamma_{\rm sonic}^1}$
and $\varphi_1<\varphi_2$ on $\overline{\Gamma_{\rm sonic}^2}$.
Then Problem \ref{problem-1} can be further reformulated into the following
free boundary problem:

\smallskip
\begin{problem}[Free Boundary Problem]\label{fbp-a}
For $\theta_{\rm w}\in (0, \theta_{\rm w}^{\rm d})$,
find a free boundary {\rm (}curved shock{\rm )} $\shock$ and a function $\vphi$ defined in domain $\Omega$,
as shown in Figs. {\rm \ref{fig:global-structure-1}}--{\rm \ref{fig:global-structure-2}},
such that $\vphi$ satisfies
\begin{itemize}
\item[\rm (i)]
Eq. \eqref{2-1} in $\Omega$,
\item[\rm (ii)]
$\vphi=\ivphi$ and $\rho \D\vphi\cdot\nnu_{\rm s}=\rho_0 \D\ivphi\cdot\nnu_{\rm s}$ {on} $\shock$,
\item[\rm (iii)]
$\vphi=\hat{\vphi}$ and $\D\vphi=\D\hat{\vphi}$ {on} $\Gamma_{\rm sonic}^1\cup\Gamma_{\rm sonic}^2\,\,$
when $\theta_{\rm w}\in (0, \theta_{\rm w}^{\rm s})$
and on $\Gamma_{\rm sonic}^2\cup \{O\}$ when $\theta_{\rm w}\in [\theta_{\rm w}^{\rm s}, \theta_{\rm w}^{\rm d})$
for $\hat{\vphi}:=\max(\vphi_1, \vphi_2)$,
\item[\rm (iv)]
$\D\vphi\cdot \nnu_{\rm w}=0$ {on} $\Wedge$,
\end{itemize}
where $\nnu_{\rm s}$ and $\nnu_{\rm w}$ are
unit normals to
$\shock$ and $\Wedge$ pointing to the interior
of $\Omega$, respectively.
\end{problem}

It can be shown that $\varphi_1>\varphi_2$ on $\Gamma_{\rm sonic}^1$, and the opposite inequality holds
on $\Gamma_{\rm sonic}^2$. This justifies the requirements in {\rm Problem \ref{fbp-a}(iii)} above.
The conditions in
{\rm Problem \ref{fbp-a}}{\rm (ii)}--{\rm (iii)}
are the Rankine-Hugoniot conditions
\eqref{1-h}--\eqref{1-i} on $\Shock$ and
$\Gamma_{\rm sonic}^1\cup\Gamma_{\rm sonic}^2$ or $\Gamma_{\rm sonic}^2\cup \{O\}$,
respectively.

\subsection{Global Solutions of Riemann Problem III:\\
Free Boundary Problem, Problem \ref{fbp-a}}

To solve Riemann Problem III, it suffices to solve the free boundary problem, Problem \ref{fbp-a},
for all the wedge angles
$\theta_{\rm w}\in (0, \theta_{\rm w}^{\rm d})$.
To obtain a global solution from $\vphi$ that is a solution of {Problem \ref{fbp-a}}
such that $\shock$ is a $C^1$--curve up to its endpoints and $\vphi\in C^1(\overline\Omega)$,
we consider two cases:

\newcommand{\SzSeg}{{{S}_{0,{\rm seg}}}}
\newcommand{\SoSeg}{{{S}_{1,{\rm seg}}}}
For the supersonic case $\theta_{\rm w}\in (0, \theta_{\rm w}^{\rm s})$,
we divide domain $\Lambda$ into four separate domains; see Fig. \ref{fig:global-structure-1}.
Denote by $\SzSeg$ the line segment $OP_1\subset S_0$, and
by $\SoSeg$ the portion (half-line) of $S_1$ with left endpoint $P_2$
so that $\SoSeg\subset\Lambda$.
Let $\Omega_{{S}}$ be the unbounded domain below
curve $\ol{\SzSeg\cup\shock\cup \SoSeg}$
and above $\Wedge$
(see Fig. \ref{fig:global-structure-1}).
In $\Omega_{{S}}$, let $\Omega_1$ be the bounded domain enclosed
by ${S}_0, \Gamma^1_{\rm sonic}$,
and $\Wedge$.
Set $\Omega_{2}:=\Omega_{{S}}\setminus \ol{\Omega_1\cup\Omega}$.
Define a function $\vphi_*$ in $\Lambda$ by
\begin{equation}\label{extsol}
\vphi_*=
\begin{cases}
\ivphi& \qquad \text{in $\Lambda\setminus \Omega_{{S}}$},\\
\vphi_1& \qquad \text{in $\Omega_1$},\\
\vphi& \qquad \text{in $\Gamma^1_{\rm sonic}\cup\Omega\cup\Gamma^2_{\rm sonic}$},\\
\vphi_2&\qquad \text{in $\Omega_2$}.
\end{cases}
\end{equation}
By {Problem \ref{fbp-a}}(ii)--(iii),
$\vphi_*$ is continuous in $\Lambda\setminus\Omega_S$
and $C^1$ in $\overline{\Omega_{{S}}}$.
In particular, $\vphi_*$ is $C^1$ across $\Gamma^1_{\rm sonic}\cup\Gamma^2_{\rm sonic}$.
Moreover, using  {Problem \ref{fbp-a}}(i)--(iii), we obtain
that  $\vphi_*$ is a global entropy solution of Eq. \eqref{2-1}
in $\Lambda$.

For the subsonic case $\theta_{\rm w}\in [\theta_{\rm w}^{\rm s}, \theta_{\rm w}^{\rm d})$,
domain $\Omega_1\cup \Gamma^1_{\rm sonic}$ in $\varphi_*$ reduces to one point $\{O\}$;
see Fig. \ref{fig:global-structure-2}.
The corresponding function $\varphi_*$ is a global entropy solution of Eq. \eqref{2-1}
in $\Lambda$.

\begin{definition}[Admissible Solutions]\label{admisSolnDef-Prandtl}
$\,$ Let $\theta_{\rm w}\in (0, \theta_{\rm w}^{\rm d})$.
A function $\varphi\in C^{0,1}(\overline\Lambda)$ is an admissible solution of
{\rm Problem \ref{fbp-a}}
if $\varphi$ is a solution of {\rm Problem \ref{fbp-a}} extended to $\Lambda$
by \eqref{extsol} and satisfies the following properties{\rm :}
\begin{enumerate}[\rm (i)]
\item\label{RegReflSol-PropPrandtl0}
The structure of solution is of the form{\rm :}
		
\smallskip
\begin{itemize}
\item
If $\theta_{\rm w}\in (0, \theta_{\rm w}^{\rm s})$,	then $\varphi$ has the configuration
shown on Fig. {\rm \ref{fig:global-structure-1}}			
such that $\Gsh$ is $C^{2}$ in its relative interior and
\begin{align*}
&\varphi\in C^{0,1}(\Lambda)\cap C^1(\Lambda\setminus (\overline{\SzSeg}\cup\overline\Gsh\cup \overline{\SoSeg})),\\
&\varphi\in  C^{1}(\overline{\Omega})\cap
				C^2(\overline\Omega\setminus(\overline{\SzSeg}\cup \overline{\SoSeg}))\cap
				C^3(\Omega).
\end{align*}

\smallskip
\item If $\theta_{\rm w}\in [ \theta_{\rm w}^{\rm s}, \theta_{\rm w}^{\rm d})$,	
then $\varphi$ has the configuration  shown on Fig. {\rm \ref{fig:global-structure-2}}
such that $\Gsh$ is $C^{2}$ in its relative interior and
\begin{align*}
&\varphi\in C^{0,1}(\Lambda)\cap C^1(\Lambda\setminus (\Gsh\cup \overline{\SoSeg})),\\
&\varphi\in  C^{1}(\overline{\Omega})\cap
				C^2(\overline\Omega\setminus(\{O\}\cup \overline{\SoSeg}))\cap
				C^3(\Omega).
\end{align*}
\end{itemize}

\smallskip		
\item\label{RegReflSol-PropPrandtl1}
Eq. \eqref{2-1} is strictly elliptic in
$\overline\Omega\setminus\,\overline{\Gso}${\rm :}
$\lvert \D\varphi\rvert<c(\lvert \D\varphi\rvert^2, \varphi)$
in $\overline\Omega\setminus\,\overline{\Gso}$.
	
\smallskip
\item\label{RegReflSol-PropPrandtl1-1}
$0<\partial_{\bn_{\rm s}}\varphi\le \partial_{\bn_{\rm s}}\varphi_0$ on $\Gsh$,
where $\bn_{\rm s}$ is the unit normal
to $\Gsh$ pointing to the interior of $\Omega$.
		
\smallskip
\item \label{RegReflSol-PropPrandtl1-1-1} The inequalities hold{\rm :}
\begin{equation}\label{phi-between-in-omega-nonSt-Prandtl}
\max\{\varphi_{1}, \varphi_{2}\}\le\varphi\le \varphi_0  \qquad \mbox{in $\Omega$}.
\end{equation}
		
\smallskip
\item\label{RegReflSol-PropPrandtl2}
The monotonicity properties hold{\rm :}
\begin{equation}\label{MonotoneProperty}
\D(\varphi_0-\varphi)\cdot \mathbf{e}_{{S}_1}\ge 0,
\quad \D(\varphi_0-\varphi)\cdot \mathbf{e}_{{S}_0}\le 0 \qquad\,\,\,\, \mbox{in $\Omega$},
\end{equation}
where $\mathbf{e}_{{S}_0}$ and $\mathbf{e}_{{S}_{1}}$ are the unit
vectors along lines ${S}_0$ and ${S}_1$
pointing to the positive $\xi_1$--direction, respectively.
\end{enumerate}
\end{definition}

The monotonicity properties in \eqref{MonotoneProperty} imply that
\begin{equation}\label{coneOfMonotPrandtlRefl-cone}
  \D(\varphi_1-\varphi)\cdot {\bf e}\le 0 \qquad\,\, \mbox{in $\overline\Omega\,\,$ for
  all ${\bf e}\in \overline{Cone(-\mathbf{e}_{{S}_1}, \mathbf{e}_{{S}_0})}$},
\end{equation}
where $Cone(-\mathbf{e}_{{S}_1}, \mathbf{e}_{{S}_0})
=\{-a\,\mathbf{e}_{{S}_1}+b \,\mathbf{e}_{{S}_0}\;: \; a, b>0\}$.
Notice that $\mathbf{e}_{{S}_0}$ and $\mathbf{e}_{{S}_1}$ are not parallel
if $\theta_{\rm w}\ne 0$.
Then we have the following theorem:

\smallskip
\begin{theorem}[Bae-Chen-Feldman \cite{BCF-14}]\label{mainPrandtlReflThm}
Let $\gamma>1$ and $u_{0}>c_0$.
For any $\theta_{\rm w}\in (0, \theta_{\rm w}^{\rm d})$, there exists a global entropy solution
$\varphi$ of {\rm Problem \ref{fbp-a}} such that the following regularity properties are satisfied for some $\alpha\in(0,1)${\rm :}

\smallskip
\begin{enumerate}[{\rm (i)}]
\item If $\theta_{\rm w}\in (0, \theta_{\rm w}^{\rm s})$,
the reflected shock $\overline{\SzSeg}\cup\Shock\cup\overline{\SoSeg}$ is $C^{2,\alpha}$--smooth,
and $\varphi\in C^{1,\alpha}(\overline\Omega)\cap C^\infty(\overline\Omega\setminus (\overline{\Gamma_{\rm sonic}^1}\cup
\overline{\Gamma_{\rm sonic}^2}))$.

\smallskip
\item If  $\theta_{\rm w}\in [\theta_{\rm w}^{\rm s}, \theta_{\rm w}^{\rm d})$,
the reflected shock $\overline{\Shock}\cup\overline{\SoSeg}$ is $C^{1,\alpha}$ near $O$
and $C^{2,\alpha}$ away from $O$, and
$\varphi\in C^{1,\alpha}(\overline\Omega)\cap C^\infty(\overline\Omega\setminus (\{O\}\cup\overline{\Gamma_{\rm sonic}^2}))$.
\end{enumerate}

\smallskip
Moreover, in both cases, $\varphi$ is $C^{1,1}$ across the sonic arcs, and
$\Shock$ is $C^\infty$ in its relative interior.
Furthermore, $\varphi$ is an admissible solution in the sense of {\rm Definition \ref{admisSolnDef-Prandtl}},
so $\varphi$ satisfies the additional properties listed in
{\rm Definition \ref{admisSolnDef-Prandtl}}.
\end{theorem}

To achieve this, for any small $\delta>0$,
the required uniform estimates of admissible solutions with wedge angles
$\theta_{\rm w} \in [0, \theta_{\rm w}^{\rm d}-\delta]$
are first obtained.
Using these estimates,
the Leray-Schauder degree theory
can be applied to obtain the existence
in the class of admissible solutions for each $\theta_{\rm w} \in [0, \theta_{\rm w}^{\rm d}-\delta]$,
starting from the unique normal solution for $\theta_{\rm w}=0$.
Since $\delta>0$ is arbitrary,
the existence of a global entropy solution for any $\theta_{\rm w}\in(0,\theta_{\rm w}^{\rm d})$
can be established.
More details can be found in Bae-Chen-Feldman \cite{BCF-14}; see also
Chen-Feldman \cite{CF-book2018} and related references cited therein.

Recently, we have also established the convexity of transonic shocks for the Prandtl-Meyer reflection configurations.

\begin{theorem}[Chen-Feldman-Xiang \cite{ChenFeldmanXiang}]\label{thm:con-Prandtl}
If a solution of the Prandtl-Meyer problem is admissible in the sense of
{\rm Definition \ref{admisSolnDef-Prandtl}},
then its domain $\Omega$ is convex, and the shock curve $\Gsh$ is a strictly convex graph.
That is, $\Gsh$ is uniformly convex on any closed subset of its relative interior.
Moreover, for the solution of {\rm Problem \ref{fbp-a}} extended to $\Lambda$
by \eqref{extsol} $($with the appropriate modification for the subsonic/sonic case$)$ with pseudo-potential
$\varphi\in C^{0,1}(\Lambda)$ satisfying
{\rm Definition \ref{admisSolnDef-Prandtl}}{\rm (\ref{RegReflSol-PropPrandtl0})}--\eqref{RegReflSol-PropPrandtl1-1-1},
the shock is strictly convex if and only if
{\rm Definition \ref{admisSolnDef-Prandtl}}\eqref{RegReflSol-PropPrandtl2} holds.
\end{theorem}

\smallskip
With the convexity of reflected-diffracted transonic shocks,
the uniqueness and stability of global regular shock reflection-diffraction
configurations have also been established  in
the class of {\em admissible solutions}; see Chen-Feldman-Xiang \cite{CFX-Unique} for the details.

\smallskip
The existence results in Bae-Chen-Feldman \cite{BCF-14}
indicate that the steady weak supersonic/transonic shock solutions are the asymptotic limits
of the dynamic self-similar solutions, the Prandtl-Meyer reflection configurations,
in the sense of \eqref{1-k-c} in Problem \ref{problem-1} for all $\theta_{\rm w}\in (0, \theta_{\rm w}^{\rm d})$
and all $\gamma>1$.

On the other hand, it is shown in Elling \cite{Elling} and  Bae-Chen-Feldman \cite{BCF-14} that,
for each $\gam >1$, there is no self-similar {\it strong} Prandtl-Meyer reflection configuration
for the unsteady potential flow in the class of admissible solutions.
This means that the situation for the dynamic stability of the steady oblique shocks of stronger strength
is more sensitive.

\section{Two-Dimensional Riemann Problem IV: \\ the von Neumann Problem for Shock Reflection-Diffraction
for the Euler Equations for Potential Flow}

In this section, we present some recent developments in the analysis of the fourth Riemann problem,
Riemann Problem IV -- the von Neumann problem for shock reflection-diffraction by wedges for the Euler equations for potential flow in form \eqref{1-a}--\eqref{1-c}, or \eqref{HCL-2}
with \eqref{1-b1}--\eqref{HCL-2e}.

\subsection{2-D Riemann Problem IV: The von Neumann Problem for Shock Reflection-Diffraction by Wedges}
\label{RegReflProblSect}

When a vertical planar shock perpendicular to the flow
direction and separating two uniform states (0) and (1),
with constant velocities
$\uu_0= (0, 0)$ and $\uu_1=(u_1, 0), u_1>0$,
and constant densities $\rho_0<\rho_1$
(state (0) is ahead or to the right of the shock, and state
(1) is behind the shock), hits a symmetric wedge W in \eqref{wedge-symm}
head-on at time $t = 0$,
a reflection-diffraction process takes place when $t > 0$.
Mathematically, the shock reflection-diffraction problem
is a 2-D lateral Riemann problem in domain
$\RR^2\setminus \overline{W}$.

\begin{problem}[Riemann Problem IV -- the von Neumann Problem for Shock Reflection-Diffraction by Wedges]\label{ibvp-c}
{\it
Piecewise constant initial data, consisting of state $(0)$
on $\{x_1>0\}\setminus \overline{W}$
and state $(1)$
on $\{x_1 < 0\}$ connected by a shock at $x_1=0$,
are prescribed at $t = 0$.
Seek a solution of Eq. \eqref{HCL-2} with \eqref{1-b1}--\eqref{HCL-2e}
for $t\ge 0$ subject to the initial
data and the boundary condition $\nabla_{\xx}\Phi\cdot\nnu_{\rm w}=0$ on $\partial W$.
}
\end{problem}

Similarly to Definition \ref{weakSol-def-Prob1}, we can define the notion of weak solutions of {Problem \ref{ibvp-c}},
by noting that
the boundary condition can be written as $\rho\nabla_{\xx}\Phi\cdot\nnu_{\rm w}=0$ on $\partial W$,
which is the spatial conormal condition for Eq. \eqref{HCL-2} with \eqref{1-b1}--\eqref{HCL-2e}.

The mathematical analysis of the shock reflection-diffraction by wedges was first proposed by
John von Neumann in \cite{Neumann1,Neumann2,Neumann}.
The complexity of reflection-diffraction configurations was first reported
by Ernst Mach \cite{Mach} in 1878, who observed two patterns of
reflection-diffraction configurations:
Regular reflection (two-shock configuration; see
Figs. \ref{figure: free boundary problems-1}--\ref{figure: free boundary problems-2})
and Mach reflection (three-shock/one-vortex-sheet configuration).
It has been found later that the reflection-diffraction configurations
can be much more complicated than
what Mach originally observed;
see also \cite{BD,CF-book2018,CFr,GlimmMajda,Guderley,VD}
and the references cited therein.

\subsection{Reformation of Riemann Problem IV}

{Problem \ref{ibvp-c}}
is invariant under self-similar scaling \eqref{sss}, so it also admits self-similar solutions determined by
Eq. \eqref{2-1}--\eqref{1-o}, along with the appropriate boundary conditions.
By the symmetry of the problem with respect to the $\xi_1$--axis,
we consider only the upper half-plane
$\{\xi_2>0\}$ and prescribe the boundary condition: $\varphi_{\bn}=0$
on the symmetry line $\{\xi_2=0\}$.
Then {Problem \ref{ibvp-c}} is reformulated as a boundary value problem
in the unbounded domain
$$
\Lambda:=\RR^2_+\setminus\{\xxi\,:\,\lvert \xi_2\rvert\le \xi_1 \tan\theta_{\rm w}, \xi_1>0\}
$$
in the self-similar coordinates,
where $\mR^2_+:=\mR^2\cap\{\xi_2>0\}$.
The incident shock in
the $\xxi$--coordinates is the half-line: $S_0=\{\xi=\xi_1^0\}\cap\Lambda$, where
\begin{equation}\label{locIncShock}
\xi_1^0:=\rho_1\sqrt{\frac{2(c_1^2-c_0^2)}{(\gamma-1)(\rho_1^2-\rho_0^2)}}=\frac{\rho_1u_1}{\rho_1-\rho_0},
\end{equation}
which is determined by the Rankine-Hugoniot conditions between states (0) and (1) on $S_0$.
Then {Problem \ref{ibvp-c}} for self-similar solutions becomes
the boundary value problem:
{\it Seek a
solution $\varphi$ of Eq. \eqref{2-1}--\eqref{1-o}
in the self-similar
domain $\Lambda$ with the slip boundary condition
$\D\varphi\cdot\nnu\rvert_{\partial\Lambda}=0$
and the asymptotic boundary condition at infinity{\rm :}
$$
\varphi\to\bar{\varphi}=
\begin{cases} \varphi_0 \qquad\mbox{for}\,\,\,
                         \xi_1>\xi_1^0, \xi_2>\xi_1 \tan\theta_{\rm w},\\
              \varphi_1 \qquad \mbox{for}\,\,\,
                          \xi_1<\xi_1^0, \;\xi_2>0,
\end{cases}
\qquad \mbox{when $\lvert \xxi\rvert\to \infty$,}
$$
where $\pSi_0=-\frac{1}{2}\lvert \xxi\rvert^2$ and $\pSi_1=-\frac{1}{2}\lvert \xxi\rvert^2+u_1(\xi_1-\xi^0_1)$.
}

Similarly, we can define the notion of weak solutions of the boundary value problem
by observing that
the boundary condition can be written as
$\rho \D\varphi\cdot\nnu\rvert_{\partial\Lambda}=0$,
which is the spatial conormal condition for Eq. \eqref{2-1}--\eqref{1-o}.
A weak solution is called an entropy solution if it satisfies
the entropy condition:
{\it density $\rho$ increases
in the pseudo-flow direction of $\D\vphi\rvert_{\Lambda^+\cap S}$
across any discontinuity curve $(${\it i.e.}, shock$)$}.

\begin{figure}[htp]
\begin{center}
	\begin{minipage}{0.42\textwidth}
		\centering
		\includegraphics[width=1.0\textwidth]{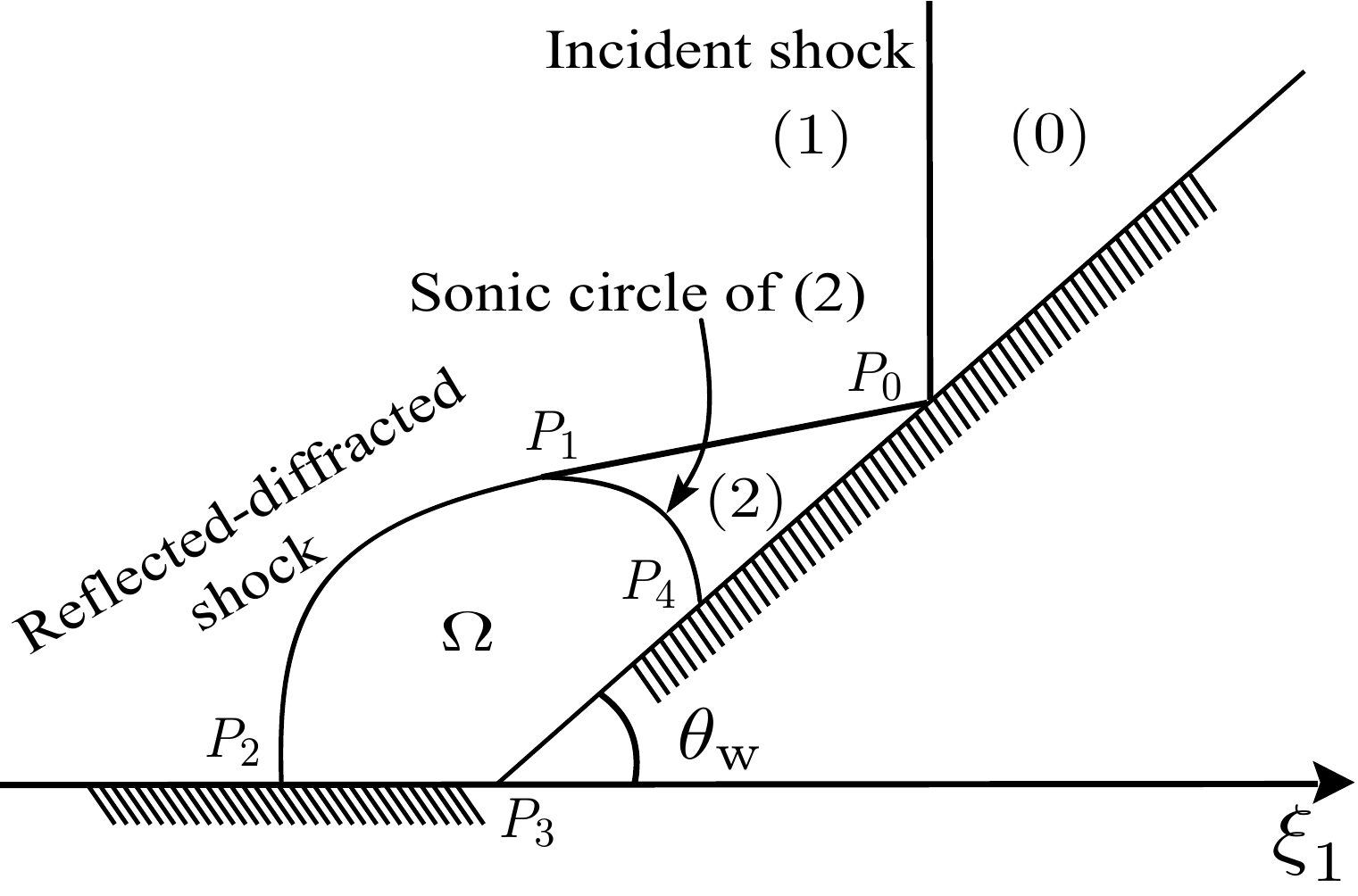}
		\caption{Supersonic regular shock reflection-diffraction configuration}
		\label{figure: free boundary problems-1}
	\end{minipage}
	\hspace{0.3in}
	\begin{minipage}{0.41\textwidth}
		\centering
		\includegraphics[width=1.0\textwidth]{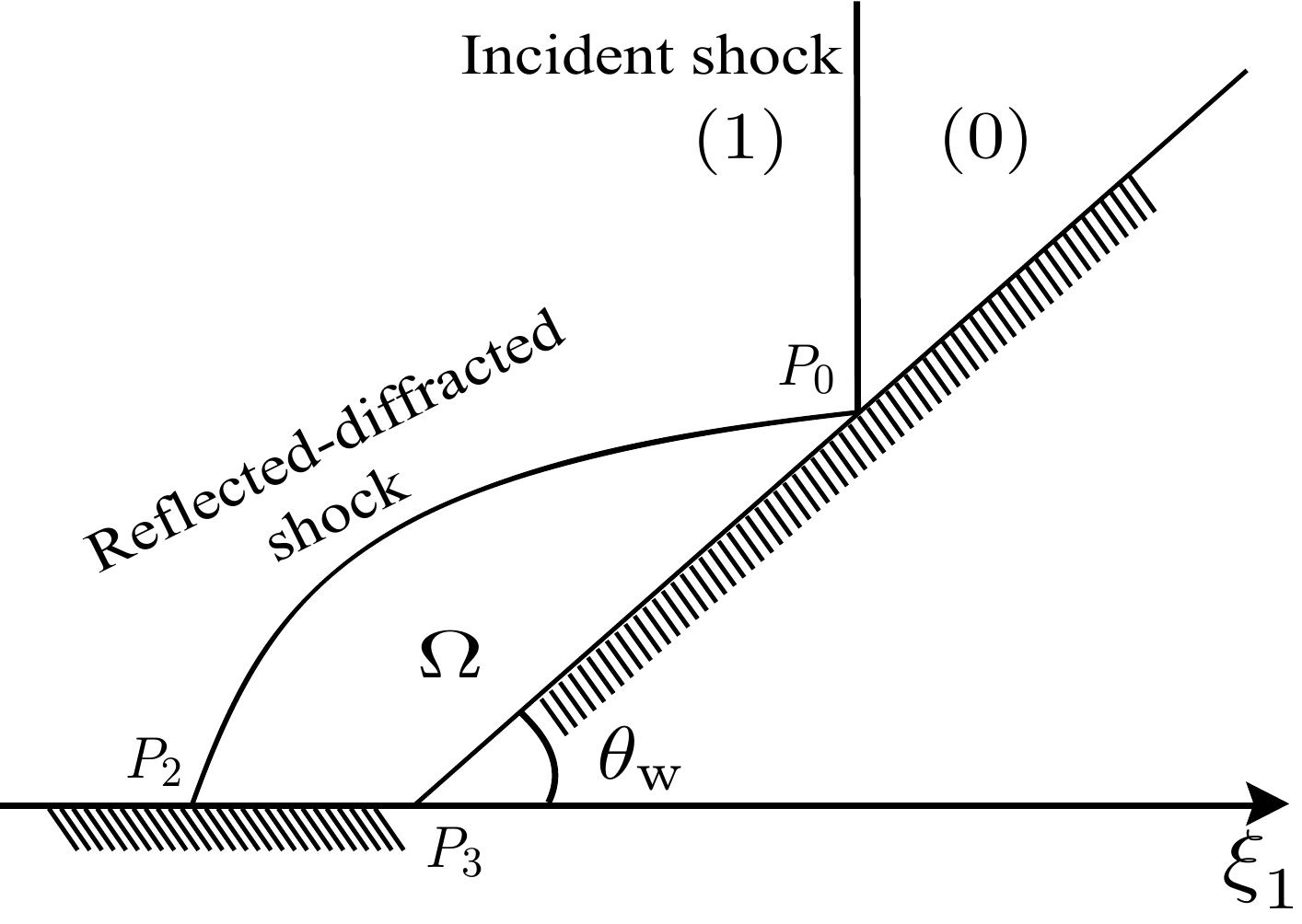}
		\caption{Subsonic regular shock reflection-diffraction configuration}
		\label{figure: free boundary problems-2}
	\end{minipage}
\end{center}
\end{figure}

\smallskip
If a solution has one of the regular shock reflection-diffraction configurations as shown
in Figs. \ref{figure: free boundary problems-1}--\ref{figure: free boundary problems-2} ({\it cf.} \cite{CF-book2018})
and its pseudo-potential
$\varphi$ is smooth in the subdomain $\Omega$
between the wedge and the reflected-diffracted shock, then
it should satisfy the slip boundary condition on the wedge and the Rankine-Hugoniot conditions
with state $(1)$ across the flat shock $S_1=\{\varphi_1=\varphi_2\}$, which passes
through point $P_0$ where the incident shock meets the wedge boundary.
Define the uniform state (2) with pseudo-potential $\varphi_2(\xxi)$ such
that
$$
\varphi_2(P_0)=\varphi(P_0), \qquad \D\varphi_2(\PtIncW)=
\lim_{P\to P_0,\; P\in \Omega} \D\varphi(P).
$$
Then the constant density $\rho_2$ of state (2) is equal to
$\rho(\lvert \D\varphi\rvert^2, \varphi)(\PtIncW)
=\rho(\lvert \D\varphi_2\rvert^2, \varphi_2)(\PtIncW)
$
via (\ref{1-o}).
It follows that
$\varphi_2$ satisfies the following three conditions at $P_0$:
\begin{equation}\label{condState2}
\D\varphi_2\cdot\bn_{\rm w}=0, \quad\varphi_2=\varphi_1,
\quad\r(\lvert \D\varphi_2\rvert^2,\varphi_2)\D\varphi_2\cdot\bn_{S_1}=
\rho_1\D\varphi_1\cdot\bn_{S_1}
\end{equation}
for $\bn_{S_1}=\frac{\D(\varphi_1-\varphi_2)}{\lvert \D(\varphi_1-\varphi_2)\rvert}$,
where  $\bn_{\rm w}$ is the outward normal to the wedge boundary.

\smallskip
State (2) can be either supersonic or subsonic at $\PtIncW$, which
determines the supersonic or subsonic type of the configurations.
The regular reflection solution in the supersonic domain is expected to
consist of the constant states separated by straight shocks ({\it cf.} \cite[Theorem 4.1]{Serre}).
Then, when state (2) is supersonic at $\PtIncW$, it can be shown that the constant state (2), extended
up to arc $\Sonic:=P_1P_4$ of the sonic circle of state (2),
as shown in Fig. \ref{figure: free boundary problems-1},
satisfies Eq. (\ref{2-1}) in the domain,
the Rankine-Hugoniot conditions \eqref{1-h}--\eqref{1-i} on the straight shock $P_0P_1$,
and the slip boundary condition: $\D\varphi_2\cdot\bn_{\rm w}=0$ on the
wedge $P_0P_4$,
and is expected to be a part of the configuration.
Then the supersonic regular shock reflection-diffraction configuration
on Fig. \ref{figure: free boundary problems-1}
consists of three uniform states (0), (1), (2), and a non-uniform state in domain
$\Omega=P_1P_2P_3P_4$,
where Eq. \eqref{2-1} is elliptic.
The elliptic domain $\Omega$ is separated from the hyperbolic domain $\PtIncW\PtUpL\PtUpR$
of state (2) by the sonic arc $\Sonic$,
on which the ellipticity in $\Omega$ degenerates.
The subsonic regular shock reflection-diffraction configuration
as shown in Fig. \ref{figure: free boundary problems-2}
consists of two uniform states (0) and (1), and a non-uniform
state in domain
$\Omega=P_0P_2P_3$,
where  Eq. \eqref{2-1}
is elliptic, and
$\varphi_{\lvert \Omega}(\PtIncW)=\varphi_2(\PtIncW)$ and
$\D(\varphi_{\rvert\Omega})(\PtIncW)=\D\varphi_2(\PtIncW)$.

For the supersonic case
in Fig. \ref{figure: free boundary problems-1},
we also use $\Shock$, $\Wedge$, and $\Symm$ for the curved part of $P_1P_2$,
the wedge boundary $P_3P_4$, and the symmetry line segment $P_2P_3$, respectively.
For the subsonic case
in Fig. \ref{figure: free boundary problems-2},
$\Shock$, $\Wedge$, and $\Symm$ denote $P_0P_2$, $P_0P_3$, and $P_2P_3$, respectively.
We unify the notations with the supersonic case by introducing
points $\PtUpL$ and $\PtUpR$ for the subsonic case as
\begin{equation}\label{P1-P4-P0-subs-Ch2}
\PtUpL\defd\PtIncW, \quad \PtUpR\defd\PtIncW, \quad \overline\Sonic\defd\{\PtIncW\}.
\end{equation}

\begin{figure}
 \centering
\includegraphics[height=42mm]{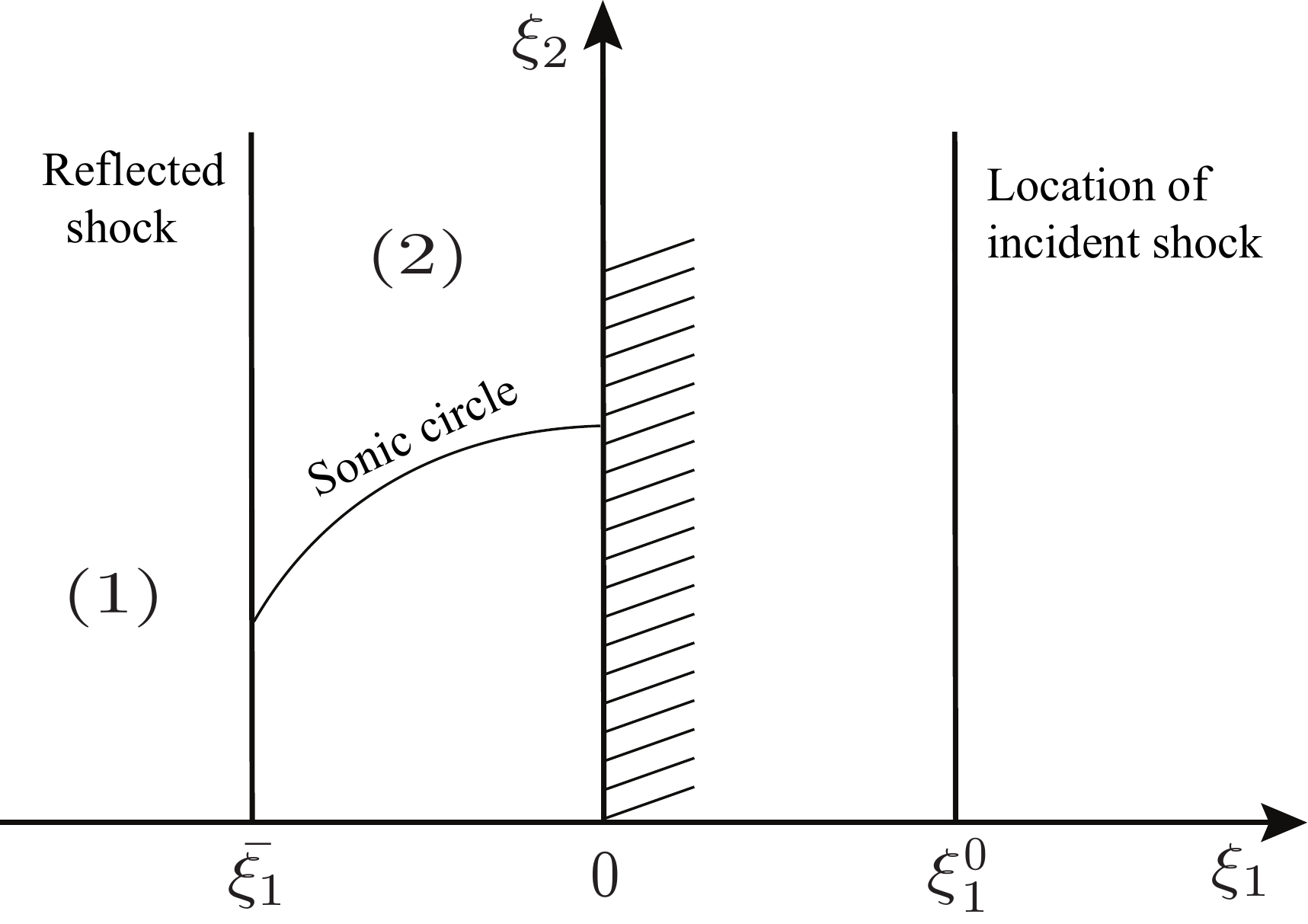}
\caption[]{Normal reflection configuration ({\it cf}. \cite{CF-book2018})}
\label{NormReflFigure}
\end{figure}

The corresponding solution for $\theta_{\rm w}=\frac\pi 2$ is called the {\em normal reflection}.
In this case, the incident shock normally reflects from the flat wall so that the reflected shock
is also a plane  $\{\xi_1=\bar{\xi}_1\}$, where $\bar{\xi}_1<0$; see Fig. \ref{NormReflFigure}.

As indicated above, a necessary condition for the existence
of a regular reflection solution is the existence of the uniform state (2) with pseudo-potential $\varphi_2$ determined by
the system of algebraic equations \eqref{condState2}
for constants $(u_2, v_2, \rho_2)$
of state (2)
across the flat shock $S_1=\{\varphi_1=\varphi_2\}$ separating it from state (1) and satisfying
the entropy conditions $\rho_2>\rho_1$.
For any fixed densities $0<\rho_0<\rho_1$ of states (0) and (1),
it can be shown that there exist a sonic angle $\theta_{\rm w}^{\rm s}$ and a detachment angle
$\theta_{\rm w}^{\rm d}$
satisfying
$$
0<\theta_{\rm w}^{\rm d}<\theta_{\rm w}^{\rm s}<\frac{\pi}{2}
$$
such that the algebraic system \eqref{condState2} has two solutions for
each $\theta_{\rm w}\in (\theta_{\rm w}^{\rm d}, \frac{\pi}{2})$ which become equal when
$\theta_{\rm w}=\theta_{\rm w}^{\rm d}$.
Thus, for each $\theta_{\rm w}\in (\theta_{\rm w}^{\rm d}, \frac{\pi}{2})$, there exist two
states (2), weak versus strong, with densities  $\rho_2^{\rm weak}<\rho_2^{\rm strong}$.
The weak state (2) is supersonic at the reflection point $\PtIncW(\theta_{\rm w})$ for
$\theta_{\rm w}\in (\theta_{\rm w}^{\rm s}, \frac{\pi}{2})$,
sonic for $\theta_{\rm w}=\theta_{\rm w}^{\rm s}$,
and subsonic for $\theta_{\rm w}\in (\theta_{\rm w}^{\rm d}, \hat\theta_{\rm w}^{\rm s})$
for some $\hat\theta_{\rm w}^{\rm s}\in(\theta_{\rm w}^{\rm d}, \theta_{\rm w}^{\rm s}]$.
The strong state (2) is subsonic at $\PtIncW(\theta_{\rm w})$ for all
$\theta_{\rm w}\in (\theta_{\rm w}^{\rm d}, \frac{\pi}{2})$.

To determine which of the two states (2) for
$\theta_{\rm w}\in (\theta_{\rm w}^{\rm d}, \frac{\pi}{2})$,
weak or strong, is physical for the local theory,
it was conjectured that the strong shock reflection-diffraction configuration
would be non-physical;
indeed, it is shown in Chen-Feldman \cite{ChenFeldman,CF-book2018}
that the weak shock reflection-diffraction configuration tends to the unique normal
reflection in Fig. \ref{NormReflFigure}, but the strong one does not,
when $\theta_{\rm w} \to \frac{\pi}{2}$.
The entropy condition and the definition of weak state (2) imply that
$0<\rho_1<\rho_2^{\rm weak}$.
With the weak state (2), the following conjectures were
proposed (see von Neumann \cite{Neumann1,Neumann2}):

\medskip
{\bf The Sonic Conjecture}:
{\em There exists a supersonic regular shock reflection-diffraction
configuration when
$\theta_{\rm w}\in (\theta_{\rm w}^{\rm s}, \frac{\pi}{2})$
for $\theta_{\rm w}^{\rm s}>\theta_{\rm w}^{\rm d}$.
That is,
the supersonicity of the weak state {\rm (2)} implies the existence
of a supersonic regular reflection
solution, as shown in Fig. {\rm \ref{figure: free boundary problems-1}.}}

\medskip
{\bf The Detachment Conjecture}:
{\em There exists a regular shock reflection-diffraction configuration for
any wedge angle $\theta_{\rm w}\in (\theta_{\rm w}^{\rm d}, \frac{\pi}{2})$.
That is, the existence of state {\rm (2)} implies the existence
of a regular reflection solution,
as shown in Figs. {\rm \ref{figure: free boundary problems-1}--\ref{figure: free boundary problems-2}}.
}

In other words, the von Neumann detachment conjecture above is that the global regular shock reflection-diffraction
configuration is possible whenever the local regular reflection at the reflection
point is possible.

From now on, for the given wedge angle
$\theta_{\rm w}\in (\theta_{\rm w}^{\rm d}, \frac{\pi}{2})$,
state (2) represents the unique weak state (2)
and $\varphi_2$ is its pseudo-potential.
State (2) is obtained from the algebraic conditions \eqref{condState2}
which determines
line $S_1$ and the sonic arc $\Sonic$ when state (2) is supersonic at $P_0$,
and the slope of $\Shock$ at $P_0$
when state (2) is subsonic at $P_0$.
Thus, the unknowns are both domain $\Omega$
and pseudo-potential $\varphi$ in $\Omega$,
as shown in Figs. {\rm \ref{figure: free boundary problems-1}}--{\rm \ref{figure: free boundary problems-2}}.
Then, from \eqref{1-h}--\eqref{1-i}, in order to construct a solution of Problem 6.1
for the supersonic or subsonic regular shock reflection-diffraction configuration, it suffices
to solve the following problem:

\begin{problem}[Free Boundary Problem]\label{fbp-c}
{\it  For $\theta_{\rm w}\in (\theta_{\rm w}^{\rm d}, \frac{\pi}{2})$,
find a free boundary $($curved reflected shock$)$ $\Shock \subset \Lambda\cap \{\cxi<\xi_{1\PtUpL}\}$
and a function $\varphi$ defined in domain
$\Omega$ as shown in Figs. {\rm \ref{figure: free boundary problems-1}}--{\rm \ref{figure: free boundary problems-2}}
such that
\begin{itemize}
\item[\rm (i)]
Eq. \eqref{2-1} is satisfied in $\Omega$ and is strictly elliptic for $\varphi$ in $\overline\Omega\setminus\overline\Sonic${\rm ,}
\item[\rm (ii)]
$\vphi=\vphi_1$ and $\rho \D\vphi\cdot\nnu_{\rm s}=\rho_1 \D\vphi_1\cdot\nnu_{\rm s}$ {on} the free boundary $\Shock${\rm ,}
\item[\rm (iii)]
$\vphi=\vphi_2$ and $\D\vphi=\D\vphi_2$ {on} $P_1P_4$
in the supersonic case as shown in Fig. {\rm \ref{figure: free boundary problems-1}}
 and at $P_0$ in the subsonic case as shown in Fig. {\rm \ref{figure: free boundary problems-1}}{\rm ,}
\item[\rm (iv)]
$\D\vphi\cdot\nnu_{\rm w}=0$ {on} $\Wedge$, and
$\D\vphi\cdot\nnu_{\rm sym}=0$ {on} $\Symm$,
\end{itemize}
where $\nnu_{\rm s}$, $\nnu_{\rm w}$, and $\nnu_{\rm sym}$ are the interior unit normals to $\Omega$
on $\shock$, $\Wedge$, and $\Symm$, respectively.
}
\end{problem}

The conditions in  {Problem \ref{fbp-c}}(ii)
are the Rankine-Hugoniot conditions \eqref{1-h}--\eqref{1-i} on $\Shock$ between
$\varphi_{\rvert\Omega}$ and $\varphi_1$.
Since $\Shock$ is a free boundary and Eq. \eqref{2-1} is strictly elliptic for $\varphi$
in $\overline\Omega\setminus\overline\Sonic$,
then two conditions (the Dirichlet and oblique derivative conditions) on $\Shock$ are consistent with
one-phase free boundary problems for nonlinear elliptic PDEs of second order.

A careful asymptotic analysis has been made for serval reflection-diffraction
configurations; see
\cite{GlimmMajda,KB,HK,Harabetian,Morawetz2}
and the references
cited therein.
Large or small scale numerical simulations have also been performed;
{\it cf.} \cite{BD,GlimmMajda,WC} and the references cited therein.
However, most of the fundamental issues for the shock reflection-diffraction
phenomena have not been understood, especially the global structures and the
transition between the different patterns of shock reflection-diffraction configurations.
This is partially because physical/numerical experiments are
hampered by many difficulties and have not yielded clear transition criteria
between the different patterns.
In particular, some different patterns occur when
the wedge angles are only fractions of a degree apart, a resolution even by
sophisticated
experiments and numerical simulations has been unable to reach ({\it cf.} \cite{BD,LD}).
Therefore, the necessary approach to understand
fully the shock reflection-diffraction phenomena, especially the transition criteria,
is via rigorous mathematical analysis.

\subsection{Global Solutions of Riemann Problem IV: \\ Free Boundary Problem, Problem  \ref{fbp-c}}

If $\varphi$ is a solution of Problem \ref{fbp-c}, define its extension
from
$\Omega$ to $\Lambda$ by
\begin{equation}\label{phi-states-0-1-2-MainThm}
\varphi=\begin{cases}
\, \varphi_0 \qquad\, \mbox{for}\,\, \xi_1>\xi_1^0 \mbox{ and } \xi_2>\xi_1\tan\theta_{\rm w},\\
\, \varphi_1 \qquad\, \mbox{for}\,\, \xi_1<\xi_1^0
  \mbox{ and above curve} \,\, P_0\PtUpL\PtLwL,\\
\, \varphi_2 \qquad\, \mbox{in domain}\,\, P_0\PtUpL\PtUpR,
\end{cases}
\end{equation}
where we have used the notational convention {\rm (\ref{P1-P4-P0-subs-Ch2})}
for the subsonic reflection case, in which
domain $P_0\PtUpL\PtUpR$ is one point
and curve $P_0\PtUpL\PtLwL$ is $P_0\PtLwL$;
see Figs. {\rm \ref{figure: free boundary problems-1}}--{\rm \ref{figure: free boundary problems-2}}.
Also,
the extension by (\ref{phi-states-0-1-2-MainThm})
is well-defined because of the requirement that $\Shock \subset \Lambda\cap \{\cxi<\xi_{1\PtUpL}\}$
in {Problem \ref{fbp-c}}.

In the supersonic case,
the conditions in {Problem \ref{fbp-c}}(iii) are the Rankine-Hugoniot conditions on $\Sonic$ between
$\varphi_{\rvert\Omega}$ and $\varphi_2$. Indeed,
since state (2) is sonic on $\Sonic$,  it follows from
\eqref{1-h}--\eqref{1-i} that no gradient jump occurs on $\Sonic$.
Then, if $\varphi$ is a solution of {Problem \ref{fbp-c}},
its extension by \eqref{phi-states-0-1-2-MainThm} is a global entropy solution in the self-similar coordinates.

Since $\Sonic$ is not a free boundary, it is not possible in general to prescribe
two conditions given in  {Problem \ref{fbp-c}}(iii) on $\Sonic$ for a second-order elliptic PDE.
In the iteration problem,
we prescribe the condition: $\varphi=\varphi_2$ on $\Sonic$, and then prove that
$\D\varphi=\D\varphi_2$ on $\Sonic$ by exploiting the elliptic degeneracy on $\Sonic$.

\medskip
The key obstacle to establish the existence of regular shock reflection-diffraction configurations
as conjectured by von Neumann \cite{Neumann1,Neumann2} is an additional possibility that,
for some wedge angle $\theta_{\rm w}^{\rm a}\in (\theta_{\rm w}^{\rm d}, \frac{\pi}2)$, shock
$\PtIncW\PtLwL$ may attach to the wedge vertex $\PtLwR$, as observed
by experimental results ({\it cf.} \cite[Fig. 238]{VD}).
To describe the conditions of such a possible attachment, we note that
$$
u_1=(\rho_1-\rho_0)
\sqrt{\frac{2(\rho_1^{\gamma-1}-\rho_0^{\gamma-1})}{\rho_1^2-\rho_0^2}}>0,\qquad  \rho_1>\rho_0,
\,\,\qquad c_1=\rho_1^{\frac{\gamma-1}2}.
$$
Then it follows from the explicit expressions above that, for each $\rho_0$,
there exists $\rho^{\rm c}>\rho_0$ such that
\begin{eqnarray*}
u_1\le c_1 \quad \mbox{if $\rho_1\in (\rho_0, \rho^{\rm c}]$}; \,\, \qquad
u_1>c_1 \quad \mbox{if $\rho_1\in (\rho^{\rm c}, \infty)$}.
\end{eqnarray*}

\smallskip
If $u_1\le c_1$, we can rule out the solution with a shock attached to the wedge vertex.
This is based on the fact that, if $u_1\le c_1$, then  the wedge vertex
$P_3=(0,0)$ lies within the sonic circle $\overline{B_{c_1}(\mathbf{u}_1)}$ of state (1),
and $\Shock$ does not intersect $\overline{B_{c_1}(\mathbf{u}_1)}$, as we show below.
If $u_1> c_1$, there would be a possibility that
the reflected shock could be attached to the wedge vertex,
as the experiments show $(${\it e.g.}, \cite[Fig. 238]{VD}$)$.

To solve the free boundary problem (Problem 6.2) involving transonic shocks for
all the wedge angles $\theta_{\rm w}\in  (\theta_{\rm w}^{\rm d},\frac{\pi}2)$,
we define the following admissible solutions.

\begin{definition}\label{admisSolnDef}
Let $\theta_{\rm w}\in (\theta_{\rm w}^{\rm d},\frac{\pi}2)$.
A function $\varphi\in C^{0,1}(\overline\Lambda)$ is an admissible solution of the regular reflection
problem if $\varphi$ is a solution of {\rm Problem \ref{fbp-c}} extended to $\Lambda$
by \eqref{phi-states-0-1-2-MainThm}
and satisfies the following properties{\rm :}

\smallskip
\begin{enumerate}[\rm (i)]
\item\label{RegReflSol-Prop0}
The structure of solution{\rm :}
		
\smallskip
\begin{itemize}
\item If $\lvert \D\varphi_2(P_0)\rvert>c_2$, then $\varphi$ is of the {\em supersonic} regular shock
reflection-diffraction configuration as shown on Fig. {\rm \ref{figure: free boundary problems-1}}	
and satisfies the conditions that the curved part of reflected-diffracted shock $\Gsh$ is  $C^{2}$ in its relative interior{\rm ;}
curves $\Gsh$,  $\Sonic$, $\Wedge$, and $\Symm$ do not have common points except their endpoints{\rm ;}
and
\begin{align*}
&\varphi\in C^{0,1}(\Lambda)\cap C^1(\Lambda\setminus ({S_0}\cup \overline{P_0P_1P_2})),\\
&\varphi\in C^{1}(\overline{\Omega})\cap C^{3}(\overline\Omega\setminus(\overline\Gso\cup\{P_2, P_3\})).
\end{align*}

\item If $\lvert \D\varphi_2(P_0)\rvert\le c_2$, then $\varphi$ is of the {\em subsonic} regular shock
reflection-diffraction configuration  shown on Fig. {\rm \ref{figure: free boundary problems-2}}			
and satisfies the conditions that the reflected-diffracted shock $\Gsh$ is  $C^{2}$ in its relative interior{\rm ;}
curves $\Gsh$,  $\Wedge$, and $\Symm$ do not have common points except their endpoints{\rm ;}
and
\begin{align*}
&\varphi\in C^{0,1}(\Lambda)\cap C^1(\Lambda\setminus ({S_0}\cup\overline\Gsh)),\\
&\varphi\in  C^{1}(\overline{\Omega})\cap
				C^3(\overline\Omega\setminus\{P_0, P_3\}).
\end{align*}
\end{itemize}
Moreover, in both the supersonic and subsonic cases,
the extended curve
$\Gsh^{\rm ext}:=\Gsh\cup \{P_0\}\cup\Gsh^-$ is $C^1$ in its relative interior,
where $\Gsh^-$ is the reflection of $\Gsh$ with respect to the $\xi_1$--axis.
		
\smallskip		
\item\label{RegReflSol-Prop1}
Eq. \eqref{2-1} is strictly elliptic in
		$\overline\Omega\setminus\,\overline{\Gso}${\rm :}
		$\lvert \D\varphi\rvert<c(\lvert \D\varphi\rvert^2, \varphi)$
		 in $\overline\Omega\setminus\,\overline{\Gso}$.
		
\smallskip
\item\label{RegReflSol-Prop1-1}
$\partial_{\bn_{\rm s}}\varphi_1>\partial_{\bn_{\rm s}}\varphi>0$ on $\Gsh$,
where $\bn$ is the normal to $\Gsh$ pointing to the interior of $\Omega$.
		
\smallskip
\item \label{RegReflSol-Prop1-1-1}  Inequalities hold{\rm :}
\begin{equation}\label{phi-between-in-omega-nonSt}
\varphi_2\le\varphi\le \varphi_1 \qquad\mbox{in $\Omega$}.
\end{equation}
		
\smallskip
\item\label{RegReflSol-Prop2}
The monontonicity properties hold:
\begin{equation}\label{coneOfMonotRegRefl}
\partial_{\xi_2}(\varphi_1-\varphi)\le 0, \quad \D(\varphi_1-\varphi)\cdot \mathbf{e}_{S_1}\le 0
\qquad\,\, \mbox{in $\Omega\,\,$ for $\mathbf{e}_{S_1}=\frac{P_0P_1}{\lvert P_0P_1\rvert}$}.
\end{equation}
\end{enumerate}
\end{definition}
Notice that \eqref{coneOfMonotRegRefl} implies that
\begin{equation}\label{coneOfMonotRegRefl-cone}
  \D(\varphi_1-\varphi)\cdot \mathbf{e}\le 0 \qquad \mbox{in $\overline\Omega\,$ for
  any $\mathbf{e}\in \overline{Cone(\mathbf{e}_{\xi_2}, \mathbf{e}_{S_1})}$},
\end{equation}
where
$Cone({\bf e}_{\xi_2}, {\bf e}_{S_1})=\{a\,{\bf e}_{\xi_2}+b \,{\bf e}_{S_1}\;:\; a, b>0\}$
with ${\bf e}_{\xi_2}=(0,1)$,
and ${\bf e}_{\xi_2}$ and ${\bf e}_{S_1}$ are not parallel if $\theta_{\rm w}\ne \frac\pi 2$.
Then we establish the following theorem:

\smallskip
\begin{theorem}[Chen-Feldman \cite{ChenFeldman,CF-book2018}]\label{mainShockReflThm}
There are two cases{\rm :}

\smallskip
\begin{enumerate}
\item[\rm (i)]
If  $\rho_0$ and $\rho_1$ are such that $u_1\le c_1$, then the
supersonic/subsonic regular reflection solution
exists for each wedge angle $\theta_{\rm w}\in (\theta_{\rm w}^{\rm d}, \frac{\pi}{2})$.
That is, for each $\theta_{\rm w}\in (\theta_{\rm w}^{\rm d}, \frac{\pi}{2})$,
there exists a solution $\varphi$ of {\rm Problem \ref{fbp-c}} such that
$$
\Phi(t, {\bf x}) =t\,\varphi(\frac{\bf x}{t}) +\frac{\lvert \bf x\rvert^2}{2t}
\qquad\mbox{for}\,\, \frac{\bf x}{t}\in \Lambda,\, t>0
$$
with
$$
\rho(t, {\bf x})=\Big(\rho_0^{\gamma-1}-(\gamma-1)\big(\Phi_t
      +\frac{1}{2}\lvert \nabla_{\bf x}\Phi\rvert^2\big)\Big)^{\frac{1}{\gamma-1}}
$$
is a global weak solution of {\rm Problem \ref{ibvp-c}} in the sense of
{\rm Definition \ref{weakSol-def-Prob1}} satisfying the entropy condition{\rm ;}
that is, $\Phi(t, {\bf x})$ is an entropy solution.

\smallskip
\item[\rm (ii)]
If  $\rho_0$ and $\rho_1$ are such that $u_1> c_1$, then there exists
$\theta_{\rm w}^{\rm a}\in [\theta_{\rm w}^{\rm d}, \frac{\pi}2)$ so that
the regular reflection solution exists for
each  wedge angle $\theta_{\rm w}\in (\theta_{\rm w}^{\rm a}, \frac{\pi}{2})$,
and the solution is of the self-similar structure
described in {\rm (i)}  above.
Moreover, if $\theta_{\rm w}^{\rm a}>\theta_{\rm w}^{\rm d}$,
then, for the wedge angle $\theta_{\rm w}=\theta_{\rm w}^{\rm a}$,
there exists an {\it attached} solution, {\it i.e.},  $\varphi$ is
a solution of {\rm Problem \ref{fbp-c}} with  $\PtLwL=\PtLwR$.
\end{enumerate}
The type of regular shock reflection-diffraction
configurations $($supersonic as in Fig. {\rm \ref{figure: free boundary problems-1}}
or  subsonic as in Fig. {\rm \ref{figure: free boundary problems-2}}$)$
is determined by the type of state
{\rm (2)} at $\PtIncW${\rm :}

\smallskip
\begin{enumerate}[\rm (a)]
\item
For the supersonic and sonic reflection case,
the reflected-diffracted shock $\PtIncW\PtLwL$ is $C^{2,\alpha}$--smooth for some $\alpha\in(0,1)$ and
its curved part $P_1P_2$ is $C^\infty$ away from $P_1$.
Solution $\varphi$ is in $C^{1,\alpha}(\overline\Omega)\cap C^\infty(\Omega)$, and is
$C^{1,1}$ across the sonic arc which is optimal{\rm ;} that is, $\varphi$ is {\em not} $C^2$ across the sonic arc.
\item
For the subsonic reflection case
$($Fig. {\rm \ref{figure: free boundary problems-2}}$)$,
the reflected-diffracted shock $\PtIncW\PtLwL$
and solution $\varphi$ in $\Omega$ is in $C^{1,\alpha}$
near $P_0$ and $P_3$ for some $\alpha\in(0,1)$, and $C^\infty$ away from $\{P_0,P_3\}$.
\end{enumerate}
Moreover, the regular reflection solution tends to the unique normal
reflection $($as in Fig. {\rm \ref{NormReflFigure})}
when the wedge angle $\theta_{\rm w}$ tends to $\frac{\pi}{2}$.
In addition, for both supersonic and subsonic reflection cases,
\begin{equation}\label{phi-between-in-omega}
 \varphi_2<\varphi<\varphi_1 \qquad\mbox{in $\Omega$}.
\end{equation}
Furthermore, $\varphi$ is an admissible solution in the sense of
{\rm Definition \ref{admisSolnDef}} below, so that $\varphi$ satisfies
further properties listed in
{\rm Definition \ref{admisSolnDef}}.
\end{theorem}

\medskip
Theorem \ref{mainShockReflThm} is proved by solving {Problem \ref{fbp-c}}.
The first results on the existence of global solutions of the free boundary problem ({Problem \ref{fbp-c}})
were obtained for the wedge angles sufficiently close to $\frac \pi 2$ in Chen-Feldman \cite{ChenFeldman}.
Later, in Chen-Feldman \cite{CF-book2018}, these results
were extended up to the detachment angle
as stated in Theorem \ref{mainShockReflThm}. For this extension, the
techniques developed in  \cite{ChenFeldman}, notably the estimates
near the sonic arc, were the starting point.
More details can be found in Chen-Feldman \cite{CF-book2018};
also see \cite{ChenFeldman}.

Furthermore, in Chen-Feldman-Xiang \cite{ChenFeldmanXiang}, we established the convexity of transonic
shocks for the regular shock reflection-diffraction configurations.

\begin{theorem}[Chen-Feldman-Xiang \cite{ChenFeldmanXiang}]\label{thm:con}
If a solution of the von Neumann problem for shock reflection-diffraction is admissible in the sense of
{\rm Definition \ref{admisSolnDef}},
then its domain $\Omega$ is convex, and the shock curve $\Gsh$ is a strictly convex graph.
That is, $\Gsh$ is uniformly convex on any closed subset of its relative interior.
Moreover, for the solution of {\rm Problem \ref{fbp-c}} extended to $\Lambda$
by \eqref{phi-states-0-1-2-MainThm}, with pseudo-potential
$\varphi\in C^{0,1}(\Lambda)$ satisfying
{\rm Definition \ref{admisSolnDef}}{\rm (\ref{RegReflSol-Prop0})}--\eqref{RegReflSol-Prop1-1-1},
the shock is strictly convex if and only if
{\rm Definition \ref{admisSolnDef}}\eqref{RegReflSol-Prop2} holds.
\end{theorem}

Furthermore, with the convexity of reflected-diffracted transonic shocks,
the uniqueness and stability of global regular shock reflection-diffraction
configurations have also been established  in
the class of {\em admissible solutions}; see Chen-Feldman-Xiang \cite{CFX-Unique} for details.

\section{Concluding Remarks}

In this paper, we have presented four different 2-D Riemann problems involving transonic shocks
through several prototypes of hyperbolic
systems of conservation laws and have showed how these Riemann problems can be formulated/solved as free boundary problems with transonic
shocks as free boundaries for the corresponding nonlinear conservation laws of mixed elliptic-hyperbolic type
and related nonlinear PDEs.
In Li-Zheng \cite{LZ2009,LZ2010}, another 2-D Riemann problem including the classical problem of the expansion of a wedge
of gas into a vacuum
for the isentropic Euler equations has also been solved; also see
the recent work  by Lai-Sheng \cite{LS} and the references cited therein on further related Riemann problems.
The other types of 2-D Riemann problems are still wide open,
even for the prototypes of hyperbolic systems of conservation laws
as discussed in this paper.

For the full Euler equations \eqref{HCL-1} with \eqref{HCL-1e},
the 2-D Riemann problems involve vortex sheets and entropy waves,
in addition to shocks and rarefaction waves;
see \cite{ChangChen,CCY1,CCY2,CCY3,CF-book2018,GlimmK,KTa,LaxLiu,li1998two,SCG,zheng2012systems} and the references cited therein.
Almost all of these Riemann problems for the full Euler equations  \eqref{HCL-1} with \eqref{HCL-1e}
are still unsolved.
In addition, all the 3-D or higher-D Riemann problems, including M-D wedge problems or M-D conic body problems,
are still open; see \cite{CCX,C-Fang,C-Fang-2,CKZ21} and the references cited therein for some recent developments
for M-D steady problems.
The nonlinear methods and related techniques/approaches originally developed in \cite{ChenFeldman1,ChenFeldman,CF-book2018}
as presented above for
solving 2-D Riemann problems involving 2-D transonic shocks should be useful in the analysis of
these longstanding Riemann problems and newly emerging problems for nonlinear PDEs;
also see \cite{Chen2,CF-book2018,ChenFeldman2022} and the references cited therein.
Certainly, further new ideas, techniques, and methods still need to be developed in order
to solve these mathematically challenging and fundamentally important problems.
\backmatter

\bmhead{Acknowledgments}
This paper is dedicated to Professor Tong Zhang (Tung Chang) on the occasion of his 90th birthday, who has been one of the pioneers and main contributors
in the analysis of the 2-D Riemann problems; see for example \cite{ChangChen,CH,CCY1,CCY2,CCY3,li1998two,TZ,ZLZ,ZZ,ZZ1} and the references cited therein.
The research of
Gui-Qiang G. Chen was supported in part by
the UK Engineering and Physical Sciences Research Council Awards
EP/L015811/1, EP/V008854/1, EP/V051121/1, and the Royal Society--Wolfson Research Merit Award WM090014.

\section*{Declarations}

There is no conflict of interest.


\end{document}